\documentclass[12pt]{article}
\usepackage{caption}\newcommand\independent{\protect\mathpalette{\protect\independenT}{\perp}}
\def\independenT#1#2{\mathrel{\rlap{$#1#2$}\mkern2mu{#1#2}}}\usepackage{graphicx,latexsym,amssymb,amsmath}
\usepackage{tikz}
\definecolor{myblue}{RGB}{173, 216, 230}
\DeclareMathOperator*{\argmin}{arg\,min}
\usepackage[french,english]{babel}
\usepackage{here} 
\def\R{{\mathbb{R}}}

\def\P{{\mathbb{P}}}

\newcommand{{\convp}}{{\buildrel p\over\longrightarrow}}

\newcommand{\red}[1]{{\color{black} #1}}

\newcommand{{\Vs}}{{\cal V}}
\newcommand{{\Ps}}{{\cal P}}
\newcommand{{\Ss}}{{\cal S}}
\newcommand{{\Xs}}{{\cal X}}
\newcommand{{\Ls}}{{\cal L}}
\newcommand{{\Ns}}{{\cal N}}
\newcommand{{\Zs}}{{\cal Z}}
\newcommand{{\Fs}}{{\cal F}}

\newtheorem{Lemma}{Lemma}[section] 
\newtheorem{Theorem}{Theorem}[section] 

\newcommand{\proofend}{$\quad\Box{~}$}

\renewcommand{\theequation} {\arabic{section}.\arabic{equation}}
\renewcommand{\baselinestretch}{1.1}

\makeatletter
\renewcommand\theequation{\thesection.\arabic{equation}}
\@addtoreset{equation}{section}
\makeatother

\usepackage{natbib}
\usepackage{url}

\begin{document}

\title{\bf Nonparametric instrumental regression with right censored duration outcomes}

\author{
{\large Jad B\textsc{eyhum}}
\footnote{ORSTAT, KU Leuven. This work was undertaken at the Toulouse School of Economics (Universit\'e Toulouse Capitole). Financial support from the European Research Council (2014-2019 / ERC grant agreement No.\ 337665) is gratefully acknowledged.}\\\texttt{\small jad.beyhum@kuleuven.be}
\and
\addtocounter{footnote}{2}
{\large Jean-Pierre F\textsc{lorens}}
\footnote{Toulouse School of Economics, Universit\'e Toulouse Capitole. Jean-Pierre Florens acknowledges funding from the French National Research Agency (ANR) under the Investments for the Future program (Investissements d'Avenir, grant ANR-17-EURE-0010).}
\\\texttt{\small jean-pierre.florens@tse-fr.eu}
\and
{\large Ingrid V\textsc{{an} K{eilegom}}}
\footnote{ORSTAT, KU Leuven. Financial support from the European Research Council (2016-2021, Horizon 2020 / ERC grant agreement No.\ 694409) is gratefully acknowledged.}\\\texttt{\small ingrid.vankeilegom@kuleuven.be}
}

\date{\today}

\maketitle

\begin{abstract}
This paper analyzes the effect of a discrete treatment $Z$ on a duration $T$. The treatment is not randomly assigned. The confounding issue is treated using a discrete instrumental variable explaining the treatment and independent of the error term of the model. Our framework is nonparametric and allows for random right censoring. This specification generates a nonlinear inverse problem and the average treatment effect is derived from its solution. We provide local and global identification properties that rely on a nonlinear system of equations. We propose an estimation procedure to solve this system and derive rates of convergence and conditions under which the estimator is asymptotically normal. When censoring makes identification fail, we develop partial identification results. Our estimators exhibit good finite sample properties in simulations. We also apply our methodology to the Illinois Reemployment Bonus Experiment.
\end{abstract}

\smallskip


\bigskip

\def\baselinestretch{1.3}

\newpage
\normalsize

\setcounter{footnote}{0}
\setcounter{equation}{0}

\section{Introduction}
\setcounter{equation}{0}

The objective of this paper is to analyze treatment models when the outcome of the treatment is a duration possibly observed with censoring.  Let $Z$ be the level of the treatment, the outcome will be a duration $T$, depending on $Z$ and on a random element $U$. If the treatment is randomly assigned, the model is formalized by the conditional distribution of $T$ given $Z$.  However, in many cases the assignment mechanism of $Z$ is not independent of $U$ and then the conditional distribution will mix the effect of the treatment and of the assignment mechanism.  In the econometric literature we are faced with a usual endogeneity problem, not specific to duration models.  The distinguishing feature of this paper is that we introduce a random right censoring mechanism.  The censoring duration is only observed for censored observations. We tackle the endogeneity issue using an instrumental variable $W$ independent of $U$ and sufficiently dependent of $Z$ given $U$. The model is nonparametric and nonseparable.

The problem is simplified by considering only the case where both $Z$ and $W$ are categorical and non-dynamic.  This avoids the question of ill-posedness of the inversion and of regularity assumptions on the functional parameters.  Under some usual conditions of completeness (see \citet{CH} and \citet{FFV}), we obtain point identification of the quantile regression function of $T$ on $Z$ in a subset of the parameter space, but outside of this set we can only identify regions of the parameter space.  This type of non-unicity of the solution of an inverse problem due to a restriction of the range of the estimated operator seems fully new.  Treatment effects on the quantiles of the duration, the survival probabilities and the hazard rates can be derived from the quantile regression function. Because we make a rank invariance assumption (see \citet{CH} and \citet{wuthrich2020comparison}), the identified treatment effects concern the whole population and not solely the subset of compliers. We propose a strategy to estimate the quantile regression function of $T$ on $Z$ by solving a nonlinear inverse problem (see \citet{D} and \citet{C}). Sufficient conditions under which the proposed estimators are asymptotically normal are derived. Inference results based on a bootstrap approach are given. Our estimation procedure exhibits good finite sample properties in simulations. We apply our methodology to the Illinois reemployment bonus experiment data.

This work borrows from the nonparametric instrumental regression  literature (\citet{Darolles} and  \citet{CH}). We extend the conventional framework to allow for censoring by using a particular conditional moment equation that can be estimated under random right censoring if the support of $T$ is included in the support of $C$ given $Z$, $W$. When the last condition fails, the regression function is only partially identified (\citet{manski1990nonparametric} and \citet{manski2003partial}). The identified set is characterized by a mix of conditional moment equalities and inequalities as in \citet{andrews2013inference}. Note also that instrumental variables are not the only method developed in econometrics to address nonparametrically the endogenity issue, see for instance the control function approach (\citet{Newey}) or the $g$-calculus (\citet{IR}).  

This paper contributes to the literature on program evaluation under censoring. Some papers introduce a `matching hypothesis' (conditional randomization) by assuming that the assignment of the treatment is conditionally independent of the outcome 
given a set of observed variables (see \citet{VdBBM, VdBBCD} and \citet{sant2016program}). In contrast, \citet{AVdB, AVdB2} specify a structural model involving separability assumptions and unobserved heterogeneity terms explaining the endogeneity of the treatment. Contrarily to the present article, their approach requires that the econometrician possesses two exogenous continuous regressors affecting both the duration and the outcome. 

Similarly to ours, other works introduce an instrumental variable to solve the endogeneity issue. \citet{chernozhukov2015quantile} consider a semiparametric framework under which, in contrast with the situation of this paper, the treatment is continuous and the censoring duration is observed. They propose an estimator and give conditions under which this estimator consistently estimates the quantile regression function of the whole population. Research in the field of biostatistics (\citet{tchetgen2015instrumental}, \citet{li2015instrumental} and \citet{chan2016reader}) has studied an additive hazard model with instrumental variables. Such an additivity assumption is not needed in our framework.

Some papers focus on the case where both the treatment and the instrumental variable are binary.  \citet{frandsen2015treatment}, \citet{sant2016program}, \citet{blanco2019bounds} make a monotonicity assumption stating that there are no defiers (see \citet{angrist1996identification}).  Unlike ours, these articles are interested in local average treatment effects on the population of compliers which they are able to estimate thanks to this monotonicity assumption (see \citet{wuthrich2020comparison}). Instead, by making a rank invariance assumption, we are able to estimate directly the quantile regression function of the whole population. Another difference between our paper and \citet{frandsen2015treatment} and \citet{sant2016program} is that the latter papers estimate the counterfactual cumulative distribution functions before inverting them to obtain estimates of the quantile treatment effects. This inversion process requires regularity assumptions that we do not make. Unlike in our paper, in \citet{frandsen2015treatment} the censoring duration is observed. \citet{blanco2019bounds} consider the same binary setting but allows for selection and endogenous censoring. Using principal stratification, it develops bounds on the quantile treatment effects. Finally, in the binary treatment and binary instrumental variable setting, \citet{BR} study a semiparametric and separable model where there is full compliance in the control group.


This paper is organized as follows.  The model specification is studied in Section 2. Identification results are derived in Section 3.  Section 4 is devoted to the estimation theory. Section 5 describes our simulations and the empirical application.  All the technical details are deferred to Appendices A, B and C. We draw conclusions in Section 6.

\section{The model}\label{sec:m}

Let us consider the following model:
\begin{eqnarray} \label{model}
T = \varphi(Z,U)=\varphi_Z(U),
\end{eqnarray}
where $T \in \R_+$, $Z$ is a categorical treatment with support $\{z_1,\dots, z_L\}$, $U$ has a unit exponential distribution and $\varphi$ belongs to  $L_{\uparrow}^2(Z,U)$, the set of mappings from $\{z_1,\dots,z_L\}\times \R_+$ to $\R_+$ that are square integrable with respect to the distribution of $(Z,U)$, and strictly increasing and differentiable in their second argument. The variable $W$ is a categorical instrumental variable with support $\{w_1,\dots, w_K\}$, which is independent of $U$. We assume that the distribution ot $U$ given $(Z,W)$ is continuous. Note that our model is identical to the usual nonseparable IV model in econometrics, except that $U$ follows an exponential distribution instead of a uniform distribution, which is more natural in the context of duration models (the cumulative hazard of the duration follows such a distribution). Such a normalization of the distribution of $U$ is necessary to obtain identification. Note that knowledge of $\varphi$ is equivalent to knowledge of the counterfactual quantile function of $T$ and the two differ solely because the distribution of $U$ is not uniform on $[0,1]$. In particular, for $z\in\{z_1,\dots, z_L\}$ and $u\in\R_+$, $\varphi_z(u)$ is the $1-e^{-u}$-quantile of the potential outcome of the duration for treatment level $z$, that is of $\varphi_z(U)$.  The variable $U$ represents both ex-ante unobserved heterogeneity, that is differences among individuals which are prior to the duration spell, and ex-post shocks occuring through the duration spell. The duration is right censored by a random variable $C$ with support in $\R_+$, and we observe $Y = \min(T,C)$, $Z$ , $W$ and $\delta = I(T \le C)$.  We suppose that $C$ is independent of $T$ given $(Z,W)$.   

\red{The main object of interest is the regression function $\varphi$. Many quantities of interest can be derived from $\varphi$. First, the quantile treatment effect (QTE) of a change in treatment from $z^0$ to $z^1$ for an individual with $U=u\in \mathbb{R}^+$, 
$ \varphi(z^1,u)-\varphi(z^0,u).$
Note that this is the QTE at the $1-e^{-u}$-quantile. Then, there is the average treatment effect of a change in treatment from $z^0$ to $z^1$,
$ \mathbb{E}[\varphi(z^1,U)-\varphi(z^0,U)].$
Consider also survival probabilities $\mathbb{P}(\varphi(z,U)\ge t)$ and hazard rates $\frac{\partial\mathbb{P}(\varphi(z,U)\le t)}{\partial t}/\mathbb{P}(\varphi(z,U)\ge t)$ at $t\in\mathbb{R}_+$ for $Z=z$.}
Let us provide two different characterizations of $\varphi$ in order to further illustrate its relevance: one in terms of the conditional survival function of $T$ and one in terms of the conditional hazard rate of $T$. First, thanks to the exclusion restriction $U\independent  W$ and to the fact that $U$ has an exponential distribution, $(\varphi_{z_\ell}(u))_{\ell=1}^L$ is a solution to the following system of equations in $\theta=(\theta_\ell)_{\ell=1}^L \in \mathbb{R}^L$:
\begin{equation} \label{Sy}
\sum_{\ell=1}^L S(\theta_\ell,z_\ell|w_k)=e^{-u}\quad \text{for}\ k=1,\dots,K,\ u\in\mathbb{R}_+,
\end{equation}
where $S(t,z|w) = \mathbb{P}(T \ge t,Z=z|W=w)$.  Indeed, 
\begin{eqnarray*}
\sum_{\ell=1}^L S(\varphi_{z_\ell}(u),z_\ell|w_k)&=& \sum_{\ell=1}^L \mathbb{P}(T\ge \varphi_{z_\ell}(u),Z=z_\ell|W=w_k)\\
&=& \sum_{\ell=1}^L \mathbb{P}(\varphi_{z_\ell}(U)\ge \varphi_{z_\ell}(u),Z=z_\ell|W=w_k)\\
&=& \sum_{\ell=1}^L \mathbb{P}(U\ge u ,Z=z_\ell|W=w_k) \\
&=& \mathbb{P}(U\ge u|W=w_k) = e^{-u} .
\end{eqnarray*}
From now on, we assume that $S$ is differentiable in its first argument. Concerning the second characterization, in the next lemma, we re-express our model in terms of the conditional hazard function. The latter is more common in the context of duration models.  

\begin{Lemma} \label{hazard}
Suppose $T=\varphi(Z,U)$, $U$ and $W$ are independent, $U \sim \mbox{Exp}(1)$, and the density $f(\cdot|z,w)$ of T given $Z=z,W=w$ exists.  Then,
$$ \sum_{\ell=1}^L \int_0^{\varphi_{z_\ell}(u)} h(s|z_\ell,w) p(z_\ell|T\ge s,w) \, ds  = u, $$
where $h(t|z,w) = f(t|z,w) / S(t|z,w)$ is the hazard function of $T$ given $Z=z,W=w$ (with $S$ being the corresponding survival function), and $p(z|T\ge t,w) = \mathbb{P}(Z=z|T \ge t,W=w)$.
\end{Lemma}
The proof is given in Appendix \ref{sec.A}. To derive the identification results, we use the characterization \eqref{Sy}.

\section{Identification}
\label{sec.id}

\subsection{Exact versus partial identification}
\label{sec.evp}
In this section, for the sake of simplicity, we assume that the (possibly infinite) upper bound of the support of the distribution of $C$ given $Z=z,W=w$ does not depend on $z\in\{z_1,\dots, z_L\}$ and $w\in\{w_1,\dots, w_K\}$ and we denote this upper bound by $c_0$. Such a case arises for instance when $C$ and $(T,Z,W)$ are independent.  This happens, in our empirical application where all durations greater than $26$ are censored and the others are not. Because of censoring, $S$ is only identified on $[0,c_0]\times \{z_1,\dots, z_L\}\times \{w_1,\dots, w_K\}$. We introduce the following lemma that characterizes $\varphi$ in terms of observables.
\begin{Lemma}
\label{lobseq}
For $k=1,\dots,K$, let $R_{k,u}(\theta)=\sum_{\ell=1}^L S(\theta_\ell\wedge c_0,z_\ell|w_k)- e^{-u}$. The following hold:
\begin{itemize} 
\item[(i)] If $(\varphi(z_\ell,u))_{\ell=1}^L\in [0,c_0)^{L}$, then 
$$ (\varphi(z_\ell,u))_{\ell=1}^L\in \Big\{\theta \in [0,c_0)^L \Big| R_{k,u}(\theta) = 0 \mbox { for all } k=1,\ldots,K \Big\};$$
\item[(ii)] If $(\varphi(z_\ell,u))_{\ell=1}^L\notin [0,c_0)^{L}$
, then 
$$ (\varphi(z_\ell,u))_{\ell=1}^L\in\Big\{\theta \in \mathbb{R}_+^L \Big| \max_{\ell=1}^L \theta_\ell \ge c_0, \min_{k=1}^K R_{k,u}(\theta) \ge 0 \Big\}.$$
\end{itemize}
\end{Lemma}

\noindent
{\bf Proof.}  Part (i) was proved in the previous section, thus we only prove (ii). We know that $\sum_{\ell=1}^L S(\varphi(z_\ell, u),z_\ell|w_k)= e^{-u}$. Note that $S(\cdot,z|w)$ is decreasing for any $z=z_\ell,\dots,z_L$ and $w=w_1,\dots, w_K$.  Hence, 
$$\sum_{\ell=1}^L S(\varphi(z_\ell, u)\wedge c_0,z_\ell|w_k)\ge\sum_{\ell=1}^L S(\varphi(z_\ell, u),z_\ell|w_k)= e^{-u}.$$
\hfill $\Box$\\
We make use of Lemma \ref{lobseq} to derive identification results  in the two different cases mentioned in the lemma. Let $u_0 =\argmin\limits_{\ell\in\{1,\dots, L\}}\varphi_{z_{\ell}}^{-1}(c_0)$ (if the support of $T$ is included in $[0,c_0)$, we set $u_0=\infty$). In Section \ref{sec.exid}, we discuss exact identification of $\varphi$ on $[0,u_0)$. On the interval $[u_0,\infty)$ (the empty set if $u_0=\infty$), $\varphi$ is only partially identified. We show how to use Lemma \ref{lobseq} (ii) to obtain an outer set to the identified set of $(\varphi_{z_{\ell}}(u))_{\ell=1}^L$ in Section \ref{sec.paid}.  We also discuss why the set in Lemma \ref{lobseq} (ii)  is not the identified set of $(\varphi_{z_{\ell}}(u))_{\ell=1}^L$ in general.  Finally, Section  \ref{sec.speid} obtains a smaller outer set of $(\varphi_{z_{\ell}}(u))_{\ell=1}^L$ than the one derived from Lemma \ref{lobseq} (ii) in a special case which nests our empirical application.

\subsection{Exact identification}
\label{sec.exid}
In this subsection, we discuss identification of $\varphi$ on $[0,u_0)$. For notational purposes, let us define $\mathcal{F}_{\downarrow}^{Z,W}$, the set of mappings from $\mathbb{R}_+\times \{z_1,\dots,z_L\} \times \{w_1,\dots,w_K\}$ to $\R_+$ which are continuous and decreasing in their first argument. We know that $\varphi$ belongs to the set of solutions of the equations
\begin{equation} \label{Sy2} A(\varphi,S)= 0,\end{equation}
where $A$ is the operator from 
$L_{\uparrow}^2(Z,U)\times \mathcal{F}^{Z,W}_\downarrow $ to the set of mappings from $[0,u_0)$ to $\mathbb{R}^K$ such that, for $\widetilde{\varphi} \in L_{\uparrow}^2(Z,U)$, $\widetilde{S}\in  \mathcal{F}^{Z,W}_\downarrow $ and $u\in [0,u_0)$,
$$A(\widetilde{\varphi},\widetilde{S})(u)= \Big(\sum_{\ell=1}^L \widetilde{S}(\widetilde{\varphi}_{z_\ell}(u),z_\ell|w_k)-e^{-u}\Big)_{k=1}^K.$$
If $\widetilde{S}\in  \mathcal{F}^{Z,W}_\downarrow $ is differentiable in its first argument, then, for all $\widetilde{\varphi} \in L_{\uparrow}^2(Z,U)$, we define $\Gamma(\widetilde{\varphi},\widetilde{S})$, the Fr\'echet derivative of $A$ in its first argument at the point $(\widetilde{\varphi},\widetilde{S})$.

\subsubsection{Local identification}
Under weak assumptions, it is possible to show that the system \eqref{Sy} has a unique solution in a neighborhood of $\varphi$. This is a local identification result.
The Fr\'echet derivative of $A$ in its first argument  at the point $(\varphi,S)$ is given by 
$$ \Big(\Gamma(\varphi,S) (u)\Big)_{k\ell}= - f \big(\varphi_{z_\ell}(u),z_\ell | w_k\big). $$
  Also, let $g(z\vert u,w)=\mathbb{P}(Z=z\vert U=u,W=w)$ and let $G(u)$ be the $K\times L$ matrix such that $G_{ k \ell}(u)=g(z_\ell | u,w_k)$.  We make the following assumption:
  \begin{itemize}
\item[\textbf{(L)}] For all $u\in \mathbb{R}_+$, $\mathrm{rank}(G(u))\ge L$.\end{itemize}

Note that this assumption is equivalent to the conditional completeness condition:
$$ \mathbb{E}(g(Z,U)|W=w,U=u)=0 \, \mbox{ for all } w,u \,  \Longrightarrow \, g=0  \mbox{ for all } g:\{z_1,\dots,z_L\}\times \R_+\mapsto \R. $$
The latter assumption (which implies that $K \ge L$), allows us to show the local identification of our model.

\begin{Theorem} \label{locident}
Under assumption (L) and assuming that $(\varphi_{z_\ell}(u))_{\ell=1}^L\in [0,c_0)^{L}$, the model is locally identified, in the sense that if $ \Gamma(\varphi,S)(\tilde{\varphi}-\varphi) \equiv 0$ for $\tilde{\varphi} \in L_{\uparrow}^2(Z,U)$, then $\tilde\varphi \equiv \varphi$.  
\end{Theorem}

\noindent
{\bf Proof.}  Write
\begin{eqnarray}
&& \sum_{\ell=1}^L( \tilde\varphi_{z_\ell}(u)- \varphi_{z_\ell}(u)) f \big(\varphi_{z_\ell}(u),z_\ell | w_k\big) = \sum_{\ell=1}^L \frac{\tilde\varphi_{z_\ell}(u)- \varphi_{z_\ell}(u)}{\varphi_{z_\ell}'(u)} g(u,z_\ell|w_k), \label{loc}
\end{eqnarray}
provided that $\varphi_{z_\ell}'(u) \neq 0$, where $g(u,z_\ell|w_k) = \varphi_{z_\ell}'(u) f \big(\varphi_{z_\ell}(u),z_\ell | w_k\big)$.  As $g(u,z_\ell|w_k)  = e^{-u} g(z_\ell|u,w_k)$, we have that \eqref{loc} equals zero if and only if
$$ \sum_{\ell=1}^L \frac{\tilde\varphi_{z_\ell}(u)- \varphi_{z_\ell}(u)}{\varphi_{z_\ell}'(u)} g(z_\ell|u,w_k) = 0 \Leftrightarrow G(u) \Big(\frac{\tilde\varphi_{z_\ell}(u)- \varphi_{z_\ell}(u)}{\varphi_{z_\ell}'(u)}\Big)^L_{\ell=1}=0,$$
 and hence it follows from assumption (L) that $\tilde\varphi \equiv \varphi$. \hfill $\Box$\\

It is useful to provide some intuition on assumption (L) in the case of our empirical application. 
If $Z$ is a binary treatment indicator and $W$ a binary instrument, the matrix $G(u)$ corresponds to
$$G(u)=\left(\begin{array}{cc} 
\mathbb{P}(Z=0\vert  U=u,W=0)& 1-\mathbb{P}(Z=0\vert  U=u,W=0)\\
\mathbb{P}(Z=0\vert  U=u,W=1)&1- \mathbb{P}(Z=0\vert  U=u,W=1)\end{array}\right).$$
Therefore, assumption (L) is satisfied if and only if $$\mathbb{P}(Z=0\vert  U=u,W=0)\ne \mathbb{P}(Z=0\vert  U=u,W=1)$$ for all $u\in \mathbb{R}_+$, that is the instrument changes the treatment probability for any value of \red{$U$}.  This is similar to Example 1 in \citet{D}.

\subsubsection{Global identification}
 
We continue with global identification. The aim is to find conditions under which the system \eqref{Sy} has a unique solution. Our presentation is borrowed from Appendix A of \citet{FFV}. Consider the density of $(U,Z)$ conditional on $W$. This density is perturbed in the direction of a function $\tilde\varphi$ and the amount of perturbation is characterized by the parameter $\mu>0$ :
$$ g_{\mu,\tilde\varphi}(u,z|w) = \Big[\varphi_z'(u) + \mu \tilde\varphi_z'(u)\Big] f\big(\varphi(z,u)+\mu \tilde\varphi(z,u),z|w\big). $$ 
Here, $f(t,z|w) = -\frac{d}{dt} S(t,z|w)$.  We make the following hypothesis:
\begin{itemize}
\item[\textbf{(G)}] If $\int_0^1 E_{g_{\mu,\tilde\varphi}}\big(\rho_\mu(Z,U) | U=u,W=w\big) d\mu=0$ for all $u,w$, then $\rho_\mu \equiv 0$, for any function $\rho_\mu:\{z_1,\dots, z_L\}\times\R_+\mapsto \R $.
\end{itemize}

Assumption (G) amounts to strong conditional completeness of $Z$ given $U$ and $W$ as introduced in \citet{CH} and studied in \citet{FFV}. Under this hypothesis we can now show the global identification of our model.  
\begin{Theorem} \label{globident}
Under assumption (G), $\varphi$ is globally identified on $[0,u_0)$, in the sense that if $ A(\tilde\varphi,S) \equiv A(\varphi,S)$ for $\tilde{\varphi} \in L_{\uparrow}^2(Z,U)$, then $\tilde\varphi \equiv \varphi$.  
\end{Theorem}
We refer to Appendix A of \citet{FFV} for the proof of this result in the continuous case.   The adaptation to the discrete case is immediate.

\subsection{Partial identification}
\label{sec.paid}
In this section, we consider pointwise identification at $u\in[u_0,\infty)$. In this case, by Lemma \ref{lobseq}, we know that \begin{equation}
\label{outer set}
(\varphi_{z_\ell}(u))_{\ell=1}^L\in\Big\{\theta \in \mathbb{R}_+^L \Big| \max_{\ell=1}^L \theta_\ell \ge c_0, \min_{k=1}^K R_{k,u}(\theta)\ge 0 \Big\}.
\end{equation}
This is an outer set to the identified set and Lemma \ref{lobseq} (ii) therefore constitutes a partial identification result. This outer set is not sharp. Section \ref{sec.speid} explains how to obtain a smaller outer set in the degenerate special case where $L=K=2$ and $S(t,1|0)=0$ for all $t\in\R_+$. The fact that \eqref{outer set} is not sharp is not due to this specific situation where $S(\cdot|z,w)$ is zero for a given couple $(z,w)\in \{z_1,\dots,z_L\} \times \{w_1,\dots,w_K\}$. Indeed, in the general setting where $S(\cdot,z|w)$ is nonzero for all values of $z$ and $w$, we give in the supplementary material a counterexample where \eqref{outer set} is not the identified set. This counterexample suggests that one reason why \eqref{outer set} is not a sharp set is that it does not take into account the constraints of monotonicity and differentiability of $\varphi$. Let us now discuss how to compute this outer set as a finite union of product of intervals.

\begin{figure}[ht]
\centering
\begin{tikzpicture}
\fill[color=myblue]
 (0,1) -- (0,2)
-- (0,2) -- (2,2)
-- (2,2) -- (2,0)
-- (2,0) -- (1,0)
-- (1,0) -- (1,1)
-- (1,1) -- (0,1) -- cycle ;
\draw[red] (0.5, 2) .. controls (1.15,1.15) .. (2,0.5);
\draw[->] (0,0) -- (2,0);
\draw (1,0) node[below] {$c_0$};
\draw (1,-0.1) -- (1,0.1);
\draw [->] (0,0) -- (0,2);
\draw (0,1) node[left] {$c_0$};
\draw (-0.1,1) -- (0.1,1);
\draw (1,-1) node{\tiny{$\mathbb{R}_+\times [c_0,\infty)\cup [c_0,\infty)\times\mathbb{R}_+$}};
\end{tikzpicture}
\begin{tikzpicture}
\fill[color=myblue]
 (0,1) -- (0,2)
-- (0,2) -- (0.5,2)
-- (0.5,2) -- (0.5,1)
-- cycle ;
\fill[color=myblue]
 (1,0) -- (2,0)
-- (2,0) -- (2,0.5)
-- (1,0.5) --(1,0)
-- cycle ;
\draw[red] (0, 2) .. controls (0.65,0.65) .. (2,0);
\draw[->] (0,0) -- (2,0);
\draw (1,0) node[below] {$c_0$};
\draw (1,-0.1) -- (1,0.1);
\draw [->] (0,0) -- (0,2);
\draw (0,1) node[left] {$c_0$};
\draw (-0.1,1) -- (0.1,1);
\draw (1,-1) node{\tiny{$[0,\bar{\theta}_1] \times [c_0,\infty) \cup [c_0,\infty) \times [0,\bar{\theta}_2]$}};
\end{tikzpicture}
\begin{tikzpicture}

\fill[color=myblue]
 (1,0) -- (2,0)
-- (2,0) -- (2,0.5)
-- (1,0.5) --(1,0)
-- cycle ;
\draw[red] (0, 0.85) .. controls (0.60,0.60) .. (2,0.25);
\draw[->] (0,0) -- (2,0);
\draw (1,0) node[below] {$c_0$};
\draw (1,-0.1) -- (1,0.1);
\draw [->] (0,0) -- (0,2);
\draw (0,1) node[left] {$c_0$};
\draw (-0.1,1) -- (0.1,1);
\draw (1,-1) node{\tiny{$[c_0,\infty) \times [0,\bar{\theta}_2]$}};
\end{tikzpicture}
\begin{tikzpicture}
\fill[color=myblue]
 (0,1) -- (0,2)
-- (0,2) -- (0.5,2)
-- (0.5,2) -- (0.5,1)
-- cycle ;
\draw[red] (0.25, 2) .. controls (0.60,0.60) .. (0.85,0);
\draw[->] (0,0) -- (2,0);
\draw (1,0) node[below] {$c_0$};
\draw (1,-0.1) -- (1,0.1);
\draw [->] (0,0) -- (0,2);
\draw (0,1) node[left] {$c_0$};
\draw (-0.1,1) -- (0.1,1);
\draw (1,-1) node{\tiny{$[0,\bar{\theta}_1] \times [c_0,\infty)$}};
\end{tikzpicture}
\captionsetup{labelsep=none}
\caption{. Outer set (in blue) under four different cases and the set $\{\theta \in \mathbb{R}_+^L |  \min\limits_{k=1,\dots, K} \sum_{\ell=1}^LS(\theta_\ell, z_\ell|w_k)= e^{-u} \}$ (in red).}
\label{fig:red}
\end{figure}
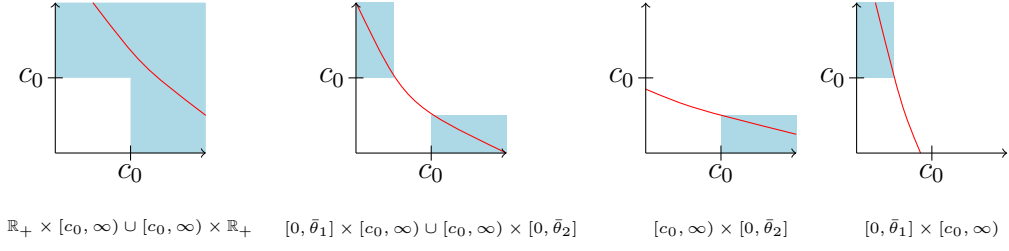

We consider first the case where $L=2$.   In that case, the set can take the following four types of shapes depending on the shape of the set $\{\theta \in \mathbb{R}_+^L | \min\limits_{k=1,\dots, K} \sum_{\ell=1}^LS(\theta_\ell, z_\ell|w_k)= e^{-u} \}$: 
\begin{itemize}
\item[] \begin{tabular}{ll}
- & the whole positive quadrant $\mathbb{R}_+^2$ with the exception of $[0,c_0)\times[0,c_0)$ \\
- & $[0,\bar\theta_1] \times [c_0,\infty) \cup [c_0,\infty) \times [0,\bar\theta_2]$ \\
- & $[c_0,\infty) \times [0,\bar\theta_2]$ \\
- & $[0,\bar\theta_1] \times [c_0,\infty)$ \\[.2cm]
\end{tabular}
\end{itemize}

\noindent
for some $0 \le \bar\theta_1,\bar\theta_2 \le c_0$. See Figure~\ref{fig:red} for a graphical description.

Let us explain how to obtain Figure~\ref{fig:red}. Take $\theta \in\R_+^2$ outside $[0,c_0)^2$. Because $S(\cdot \wedge c_0,z|w)$ is flat on $[c_0,\infty)$, $\theta$ belongs to the outer set if and only if $(\theta_1+ (c_0-\theta_1)I(\theta_1>c_0), \theta_2+ (c_0-\theta_2)I(\theta_2>c_0))$ is in the outer set. Therefore, we need to inspect the value of $\theta \rightarrow \min_{k=1}^K R_{k,u}(\theta)$ on $[0,c_0) \times \{c_0\}\cup\{c_0\}\times [0,c_0)$ to compute the outer set.  Let us start with $[0,c_0)\times \{c_0\}$. Because $S(\cdot \wedge c_0,z|w)$ is decreasing and continuous, the set of pairs $(\theta_1,c_0)$
in $[0,c_0)\times \{c_0\}$ for which $\min_{k=1}^K R_{k,u}(\theta)\ge 0$ either form a line segment $[0,\bar{\theta}_1]\times{c_0}$  where $\bar{\theta}_1\in[0,c_0]$ or the empty set. First, we consider the case where the latter set is not empty. Because $S(\cdot \wedge c_0,z|w)$ is flat on $[c_0,\infty)$, $[0,\bar{\theta}_1]\times[c_0,\infty)$ belongs to the outer set. If instead the set is empty, $[0,c_0)\times[c_0,\infty)$ does not belong to the outer set. Then, apply the same procedure on $ \{c_0\}\times [0,c_0)$. Finally, if $\min_{k=1}^K R_{k,u}((c_0,c_0)^\top)\ge 0$, $[c_0,\infty)\times [c_0,\infty)$ also belongs to the outer set.

Next, we proceed by a recursion argument and explain how the identified set can be found in dimension $L$ if it is known how to find it in dimension $L-1$.  Note that in $L$ dimensions, we can start by fixing the first coordinate $\theta_1$ to $c_0$ and specify the identified set for the remaining $L-1$ coordinates.  If the so-selected set is not empty, it should be extrapolated by letting $\theta_1$ belong to $[c_0,\infty)$.  This procedure should be repeated $L$ times by fixing each time one coordinate, and at the end the union of all obtained sets should be selected.

\subsection{A special case}
\label{sec.speid}
We conclude our discussion of identification with the special case where $L=K=2$ and $\P(Z=z_2|W=w_1)=0$. We have a triangular system of equations of the form 
\begin{eqnarray} \label{trian}
\left\{
\begin{array}{l}
S(\theta_1 ,z_1|w_1) = e^{-u} \\
S(\theta_1 ,z_1|w_2) + S(\theta_2 ,z_2|w_2) = e^{-u}.
\end{array}
\right. 
\end{eqnarray}  
This particular setting corresponds to our empirical application. In this case, if $(\varphi_{z_\ell}(u))_{\ell=1}^L\notin [0,c_0)^{L}$ and the system \eqref{Sy} has a unique solution (which for such a triangular system happens when $S$ is strictly decreasing and continuous), it is possible to obtain a smaller outer set to the identified set than $$ \Big\{\theta \in \mathbb{R}_+^L \Big| \max_{\ell=1}^L \theta_\ell \ge c_0, \min_{k=1}^K R_{k,u}(\theta) \ge 0 \Big\}.$$

We need to consider two cases.  First, if $\varphi_{z_1}(u) < c_0$, then $\theta_1=\varphi_{z_1}(u)$  is the (unique) solution of the first equation of (\ref{trian}).   This value can then be inserted in the second equation, which will only contain $\theta_2$.  If this equation has a solution, it should be $\varphi_{z_2}(u)$ and we are done.  If it does not have a solution, the identified set contains all pairs $(\varphi_{z_1}(u),\theta_2)$ with $\theta_2 \ge c_0$.  

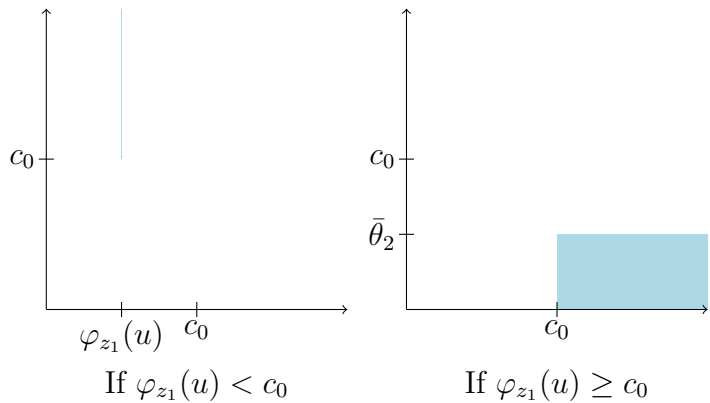
\begin{figure}[ht]
\centering
\begin{tikzpicture}

\draw[myblue] (1,2) -- (1,4) ;

\draw[->] (0,0) -- (4,0);
\draw (1,0) node[below] {$\varphi_{z_1}(u)$};
\draw (1,-0.1) -- (1,0.1);
\draw (2,0) node[below] {$c_0$};
\draw (2,-0.1) -- (2,0.1);
\draw [->] (0,0) -- (0,4);
\draw (0,2) node[left] {$c_0$};
\draw (-0.1,2) -- (0.1,2);
\draw (2,-1) node {If $\varphi_{z_1}(u)<c_0$};
\end{tikzpicture}
\begin{tikzpicture}
\fill[color=myblue]
 (2,0) -- (4,0)
-- (4,0) -- (4,1)
-- (2,1) --(2,0)
-- cycle ;
\draw[->] (0,0) -- (4,0);
\draw (2,0) node[below] {$c_0$};
\draw (2,-0.1) -- (2,0.1);
\draw [->] (0,0) -- (0,4);
\draw (0,1) node[left] {$\bar\theta_2$};
\draw (-0.1,1) -- (0.1,1);
\draw (0,2) node[left] {$c_0$};
\draw (-0.1,2) -- (0.1,2);
\draw (2,-1) node {If $\varphi_{z_1}(u)\ge c_0$};
\end{tikzpicture}
\captionsetup{labelsep=none}
\caption{. Outer set (in blue) under two different cases.}

\label{fig:red2}
\end{figure}

In the second case $\varphi_{z_1}(u) \ge c_0$.  Then, the first equation does not have a solution, and we need to solve the system of inequalities
$$ \left\{
\begin{array}{l}
S(c_0,z_1|w_1) \ge e^{-u} \\
S(c_0,z_1|w_2) + S(\theta_2 \wedge c_0,z_2|w_2) \ge e^{-u},
\end{array}
\right. $$
which leads to the identified set $[c_0,\infty) \times [0,\bar\theta_2]$, where $\bar\theta_2$ is the largest value of $\theta_2$ for which $S(\theta_2 \wedge c_0,z_2|w_2) \ge e^{-u}-S(c_0,z_1|w_2)$.  See Figure~\ref{fig:red2} for a graphical illustration of the two aforementioned cases. Note that $\theta_2$ can be equal to infinity which implies that the height of the rectangle in the second part of Figure~\ref{fig:red2} can be greater than $c_0$.

Other special cases can be considered depending on which terms in the system of equations (\ref{Sy}) are absent. However, they require a case to case analysis, which we will not further develop here.

\section{Estimation}
\label{sec.est}

\subsection{Estimation procedure}
\label{sec.estp}
We consider estimation with an i.i.d.\ sample of size $n$, $\{Y_i,Z_i,W_i,\delta_i\}_{i=1}^n$. We assume that we have an estimator $\widehat{S}$ of $S$ that satisfies properties to be specified later. Choices of $\widehat{S}$ are discussed in Section \ref{sec.choice}.
Let $\bar{U}<\infty$ be an upper bound up to which we wish to estimate $\varphi$. Let also $\bar{T}<\infty$ be an upper bound on $\max\limits_{\ell=1,\dots,L}\varphi(z_{\ell},\bar{U})$.  Let $\mathcal{K}$ be the set of mappings from $[0,\bar{U}]$ to $\mathbb{R}^K$ and $\mathcal{F}_Z^{\bar{U},\bar{T}}= \{f: \{z_1,\dots,z_L\}\times [0,\bar{U}]\mapsto [0,\bar{T}]\}$. We will make assumptions guaranteeing that $\widehat{S}$ consistently estimates $S$ on $[0,\bar{T}]$ and that the system \eqref{Sy} has a unique solution in $\mathcal{F}_Z^{\bar{U},\bar{T}}$. This implies that we are implicitly in the case where $(\varphi_{z_\ell}(u))_{\ell=1}^L$ is globally identified.

We introduce further notations. For $u\in [0,\bar{U}]$, let $V(u)$ be a positive definite $K\times K$ weighting matrix. For a $K\times K$ matrix $\bar V$ and a vector $v\in \mathbb{R}^K$, we define $\|v\|=\sqrt{v^{\top}v}$ and $\|v\|_{\bar V}=\sqrt{v^{\top}\bar Vv}$. Then, for $g\in \mathcal{K}$, we introduce $\|g\|^2=\int_{0}^{\bar{U}}\|g(u)\|^2du$, $\|g\|_{V}^2=\int_{0}^{\bar{U}}\|g(u)\|_{V(u)}^2du$. For $f\in \mathcal{F}_Z^{\bar U, \bar T}$, we use $\|f\|^2=\sup_{u\in[0,\bar{U}]}\|(f(z_{\ell},u))_{\ell=1}^L\|^2$. When we write that a random process $R$ is $O_P(a_n)$ or $o_P(a_n)$, we mean that $\|R\|=O_P(a_n)$ or $\|R\|=o_P(a_n)$, respectively. Let us now define the estimator
\begin{equation}\label{prog}  
\widehat{\varphi} \in \argmin_{\theta\in \mathcal{F}_Z^{\bar{U},\bar{T}} } \big\|A(\theta, \widehat{S})\big\|_{V}^2.
\end{equation}

\subsection{Consistency}

We introduce the following assumptions:
 \begin{itemize}
\item[\textbf{(V)}] For $u\in[0,\bar{U}]$, the eigenvalues of $V(u)$ are bounded from above and from below uniformly in $u$. The bounds are strictly positive constants.
\item[\textbf{(C)}] 
\begin{itemize}
\item[(i)]  For all $\epsilon>0$, there exists $\nu>0$ such that 
 $$\inf_{\theta\in\mathcal{F}_{Z}^{\bar{U},\bar{T}} : \|\theta-\varphi\|\ge \nu} \big\|A(\theta,S)\big\|_V\ge \epsilon ;$$
\item[(ii)] There exist $\nu,c>0$ such that for any $\theta\in \mathcal{F}_Z^{\bar{U},\bar{T}}$ satisfying $\|\theta-\varphi \| \le \nu$, 
$$ \| \theta-\varphi \| \le c \| A(\theta,S) \|;$$
\item[(iii)] There exists $(r_n)_n$, a real-valued sequence such that $r_n\to 0$ and
$$\sup_{t\in[0,\bar{U}],\ z=z_1,\dots,z_L,\ w=w_1,\dots,w_K}\big|\widehat{S}(t,z|w)-S(t,z|w)\big|=O_P(r_n).$$
\end{itemize}
\end{itemize}

Assumptions C(i) and C(ii) are conditions related to the shape of the objective function $A$ and  they ensure that there is a unique solution to the system of equations \eqref{Sy} in $[0,\bar{U}]$ and, hence, that there is a unique minimum to the program \eqref{prog} if $\widehat{S}$ estimates $S$ well enough. One can develop primitive conditions for C(ii) involving the assumption that the two first derivatives of $S$ with respect to its first argument are uniformly bounded on $[0,\bar{T}]$ (see Lemma \ref{sufficient conditions} in Appendix \ref{app.C}).  Assumption (C)(iii) implies that $\bar T$ is chosen lower than the minimum of the upper bound of the support of $C$ over all $z\in\{z_1,\dots, z_L\}$ and $w\in\{w_1,\dots, w_K\}$ such that $\mathbb{P}(Z=z|W=w)>0$. We prove the following theorem.  The proof is in Appendix A. 

\begin{Theorem} \label{C}
Under Assumptions (V) and (C), we have $\| \widehat{\varphi} -\varphi \|=O_P(r_n)$.
\end{Theorem}

This theorem shows that the rate of convergence of $\widehat{\varphi}$ is the same as the one of $\widehat{S}$ in sup-norm. Hence, solving program \eqref{prog} does not deteriorate the convergence rate. 

\subsection{Asymptotic normality}

We make the following assumption:
 \begin{itemize}
\item[\textbf{(N)}] \begin{itemize}
\item[(i)]  There exists a variance operator $\Omega$ such that $n^{1/2} A(\varphi,\widehat{S})$ converges weakly to a mean zero Gaussian process with variance operator $\Omega$;
  \item[(ii)] $\widehat{S}$ is differentiable in its first argument for $u\in[0,\bar{U}]$, the mapping $\Gamma(\varphi,S)^\top V\Gamma(\varphi,S)$ is invertible on $[0,\bar{U}]$, $ \Gamma(\varphi,\widehat{S})\xrightarrow{P}\Sigma=\Gamma(\varphi,S)$ and $ \Gamma(\widehat{\varphi},\widehat{S})\xrightarrow{P}\Sigma$;
    \item[(iii)] $A(\widehat{\varphi},\widehat{S})-A(\varphi,\widehat{S})-\Gamma(\varphi,\widehat{S})(\widehat{\varphi}-\varphi)=O_P(\| \widehat{\varphi}-\varphi \|^2);$
    \item[(iv)] $\underset{t\in[0,\bar{U}],\ z=z_1,\dots,z_L,\ w=w_1,\dots,w_K}{\sup}|\widehat{S}(t,z|w)-S(t,z|w)|=o_P(n^{-1/4})$.
\end{itemize}
\end{itemize}

This assumption leads to the following theorem, of which the proof can be found in Appendix A.  

\begin{Theorem} \label{AS}
Under Assumptions (V), (C) and (N), $\sqrt{n}(\widehat{\varphi} -\varphi )$ converges to a mean zero Gaussian process with variance operator
$[\Sigma^\top V\Sigma]^{-1} \Sigma^\top V \Omega V\Sigma[\Sigma^\top V\Sigma]^{-1} .$
\end{Theorem}
 
\subsection{Inference with bootstrap and optimality}

We propose a bootstrap procedure for inference on the functional $f((\varphi_{z_\ell}(u))_{\ell=1}^L)$ of $(\varphi_{z_\ell}(u))_{\ell=1}^L$ for $u\in\R_+$ and $f:\R_+^L\mapsto \R$ differentiable. Examples of such functionals are mentioned in Section \ref{sec:m}.
 \begin{itemize}
 \item[(1)] Draw with replacement $B$ resamples of size $n$ from the original data $\{Y_i,Z_i,W_i,\delta_i\}_{i=1}^n$;
 \item[(2)] For each resample $b$, compute $\widehat{\varphi}_b$, the value of the estimator $\widehat{\varphi}$ in the resample;
  \item[(3)] Construct a $95\%$ confidence interval for $f((\varphi_{z_\ell}(u))_{\ell=1}^L)$ using for the lower bound the $2.5\%$ percentile of $\{f((\widehat{\varphi}_b(z_\ell,u))_{\ell=1}^L)\}_{b=1}^B$ and for the upper bound the $97.5\%$ percentile of $\{f((\widehat{\varphi}_b(z_\ell,u))_{\ell=1}^L)\}_{b=1}^B$.
 \end{itemize}
\red{Lemma \ref{boot} in Appendix B proves that such a procedure works when $S$ is estimated with kernel smoothing as described in Section 4.5. The outlined procedure shows how to build pointwise confidence intervals. However, as demonstrated in Lemma \ref{boot}, a similar bootstrap approach also allows to obtain uniform confidence bands on $[0,\bar{U}]$.} 

The optimal value of $V$ is $\Sigma^{-1}$ which can be estimated using $\widehat{\Sigma}=\Gamma(\widehat{\varphi},\widehat{S})$ where $\widehat{\Sigma}$ is obtained from a first-step estimator which sets $V(u)$ equal to the identity matrix for any $u\in [0,\bar{U}]$. 

\subsection{Choices of $\widehat{\boldsymbol{S}}$}
\label{sec.choice}

In this subsection, we discuss choices for $\widehat{S}$ which are robust to random right censoring. Remark that $S(t,z\vert w)=S(t\vert z,w)p_{z,w}$ where $S(t\vert z,w)=\mathbb{P}(T\ge t\vert Z=z, W=w) $ and $p_{z,w}=\mathbb{P}(Z=z\vert W=w )$.  We define the following stochastic processes:
\begin{align*}
N_{z,w}(t)&= \sum_{i=1}^nI(Y_i\le t, Z_i=z,W_i=w,\delta_i=1); \quad & Y_{z,w}(t)&=\sum_{i=1}^nI(Y_i\ge t,Z_i=z,W_i=w);\\
Y_{z,w}&= \sum_{i=1}^nI(Z_i=z,W_i=w); & Y_{w}&=\sum_{i=1}^nI(W_i=w).
\end{align*}
We estimate $S(\cdot\vert z,w)$ using the Kaplan-Meier estimator 
\begin{equation} \label{KM} 
\widehat{S}_{KM}(t\vert z, w)=\prod_{s\le t}\Big(1-\frac{dN_{z,w}(s)}{Y_{z,w}(s)}\Big). \end{equation}
To provide an estimator $\widehat{S}$ satisfying assumption (N), one needs to smooth $\widehat{S}_{KM}$. Various techniques are available in the literature, including local polynomials and kernel smoothing. For instance, concerning the latter, if we use a kernel $K$ with a bandwidth $\epsilon$, we obtain
\begin{equation}\label{smooth} \widehat{\widetilde{S}}(t\vert z,w)=\int \widehat{S}_{KM}(t-s\epsilon\vert z, w)K(s)ds.\end{equation}
Our final estimator of $S$ is \begin{equation}\label{KM}\widehat{S}(t,z\vert w)= \widehat{\widetilde{S}}(t\vert z,w)\widehat{p}_{zw},\end{equation} where $\widehat{p}_{zw} =Y_{z,w}/Y_w$. In Appendix B, we state assumptions under which this choice of $\widehat{S}$ satisfies the conditions of Sections 4.2 and 4.3, and for which the bootstrap procedure of Section 4.4 works.

\subsection{Practical implementation} 
\label{sec.pract}
Let us now discuss how to use the outlined estimation procedure in practice. Fist, one needs to choose $\bar T$. The latter has to be chosen such that Assumption (C)(iii) holds. In appendix B, the condition on $\bar T$ guaranteeing  uniform convergence of $\widehat S$ when it is smoothed by kernel is that there exists $\xi>0$ for which $S(\bar T,z|w)/\mathbb{P}(Z=z|W=w)>\xi$ for all $z\in\{z_1,\dots, z_L\}$ and $w\in\{w_1,\dots, w_K\}$ such that $\mathbb{P}(Z=z|W=w)>0$ (condition (K)(v)).  As in Section \ref{sec.evp}, let us consider the case where the upper bound $c_0$ of the support of the distribution of $C$ given $Z=z,W=w$ does not depend on $z\in\{z_1,\dots, z_L\}$ and $w\in\{w_1,\dots, w_K\}$ and is finite. Then, if $[0,c_0)$ is strictly included in the support of $T$ (or $Y$) given $Z=z,W=w$ for all $z\in\{z_1,\dots, z_L\}$ and $w\in\{w_1,\dots, w_K\}$, the choice $\bar T=c_0$ ensures that condition (K)(v) holds. If the upper bound of the support of $C$ is $\infty$, then a simple valid choice for $\bar T$ is the minimum of the $\alpha$-quantile of $Y$ given $Z=z$ and $W=w$ over all $z\in\{z_1,\dots, z_L\}$ and $w\in\{w_1,\dots, w_K\}$ such that $\mathbb{P}(Z=z|W=w)>0$, where $\alpha\in (0,1)$. One may use $\alpha=0.95$ by convention.

In practice, the minimization program \eqref{prog} is not feasible. Instead we choose $u_1,\dots,u_M$, $M$ values at which we want to estimate $\varphi$ (for instance, $M$ i.i.d.\ replications of $U\sim \mathrm{Exp}(1)$ or a grid). For $m=1,\dots,M$, we estimate $\varphi(z,u_m)$ by $ \widehat{\varphi}(z,u_m)$ where 
\begin{equation}\label{empirical_system}
(\widehat{\varphi}(z_{\ell},u_m))_{\ell=1}^L\in \argmin_{\theta \in [0,\bar{T}]^L}  \big\|A(\theta, \widehat{S})(u_m) \big\|_{V(u_m)}^2.
\end{equation}
We obtain an estimate of $\varphi$ at the points $u_1,\dots, u_M$, which enables us to plot $\widehat{\varphi}(z,\cdot)$. 

When $\varphi$ is not identified on all its support, the estimator $(\widehat \varphi(z_\ell, u_m))_{\ell=1}^L$ will not be consistent for some values of $u$. Under the conditions of Theorem \ref{C}, by the continuous mapping theorem, $\big\|A(\widehat \varphi, \widehat{S})(u_m) \big\|_{V(u_m)}^2=o_P(1)$.  Therefore, large values of $\big\|A(\widehat \varphi, \widehat{S})(u_m) \big\|_{V(u_m)}^2$ indicate that $(\widehat \varphi(z_\ell, u_m))_{\ell=1}^L$ is not identified. As a rule of thumb, we suggest to start using partial identification results from Sections \ref{sec.paid} and \ref{sec.speid}, when $\big\|A(\widehat \varphi, \widehat{S})(u_m) \big\|_{V(u_m)}^2$ starts increasing significantly with $u$. This approach will be illustrated in Section \ref{sec.ni}.
\section{Numerical illustrations}
\label{sec.ni}

We present two applications of our approach. One corresponds to a simulated example while the second one is an experiment concerning a bonus for finding a job proposed to job seekers in Illinois. 

\subsection{Simulations}

We consider the following model.   Let $W$ follow a Bernoulli distribution with parameter 0.7. To model the dependence  between $Z$ and $(W,U)$ we set 
\begin{equation} \label{mod}
Z = I (-0.7+\varepsilon+W+0.5U\ge 0)I(W=1),
\end{equation}
where $\varepsilon\sim \mathcal{N}(0,1)$ and $U \sim Exp(1)$. This definition yields the following conditional treatment probabilities: $\mathbb{P}(Z=1|W=0)=0$ and $\mathbb{P}(Z=1|W=1)=0.76$ (average over 1,000,000 Monte Carlo replications).  Note that under this design, $W=0$ implies that $Z=0$, which is chosen to mimic the setting of our empirical application.  Next, let $T$ have an exponential distribution with hazard rate \red{$1/10$} if $Z=0$, and \red{$1/5$} if $Z= 1$, such that for $u\in\mathbb{R}_+$, $\varphi(0,u)=10u$ and $\varphi(1,u)=5u$. For $u\in\R_+$, the $1-e^{-u}$-quantile treatment effect is equal to $-5u$. Let $C=\max(15Exp(1), 10)$. Note that this implies that $C$ has a constant hazard rate of $1/15$ for durations lower than $10$. With such a censoring mechanism, around $22\%$ of the observations are censored (average of 1,000,000 replications). Finally, with $Y=\min(T,C)$ and $\delta=I(T \le C)$ we generate an i.i.d. sample $\{Y_i, Z_i,W_i, \delta_i\}_{i=1}^n$ of 10,000 observations having the same distribution as $(Y,Z,W,\delta)$.  The sample size is chosen to be close to that of the empirical example.   Note that under this data generating process, $S$ is identified on $\mathbb{R}_+$ and therefore we aim at point estimation of $\varphi$ and of the quantile treatment effects as in Section \ref{sec.est}.

$\bar{T}$ is equal to $10$, the upper bound of the support of $C$. We compute the estimators on a grid for $u$ between $0.01$ and $1.2$ with step size $0.01$. An Epanechnikov kernel is used to smooth the Kaplan-Meier estimator of the survival function. The bandwidth is chosen using the usual rule of thumb for normal densities. All results are averages over $1,000$ replications.

\begin{figure}[H]

  \centering
  \includegraphics[width=80mm]{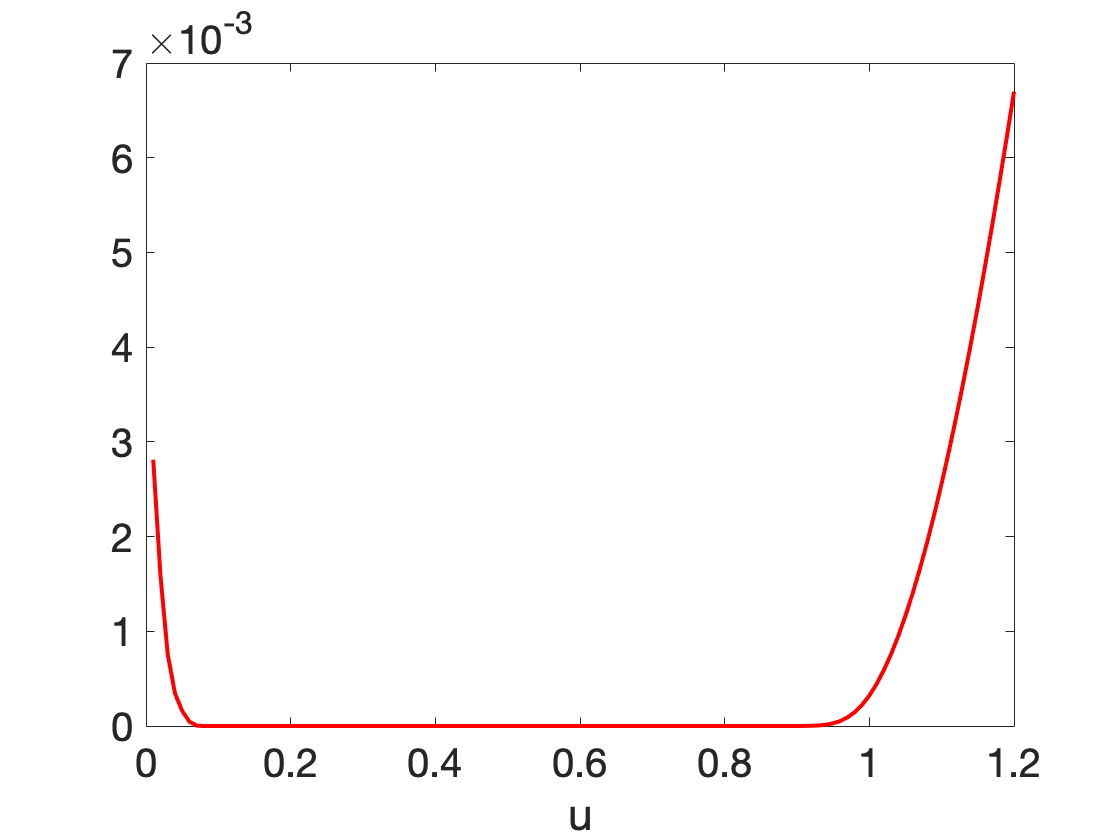}
  \caption{Value of $\big\|A(\widehat \varphi, \widehat{S})(u) \big\|_{2}^2$.}
   \label{fig:remainder_sim}

\end{figure}

\begin{figure}[H]
\begin{minipage}{0.48\textwidth}
  \centering
  \includegraphics[width=80mm]{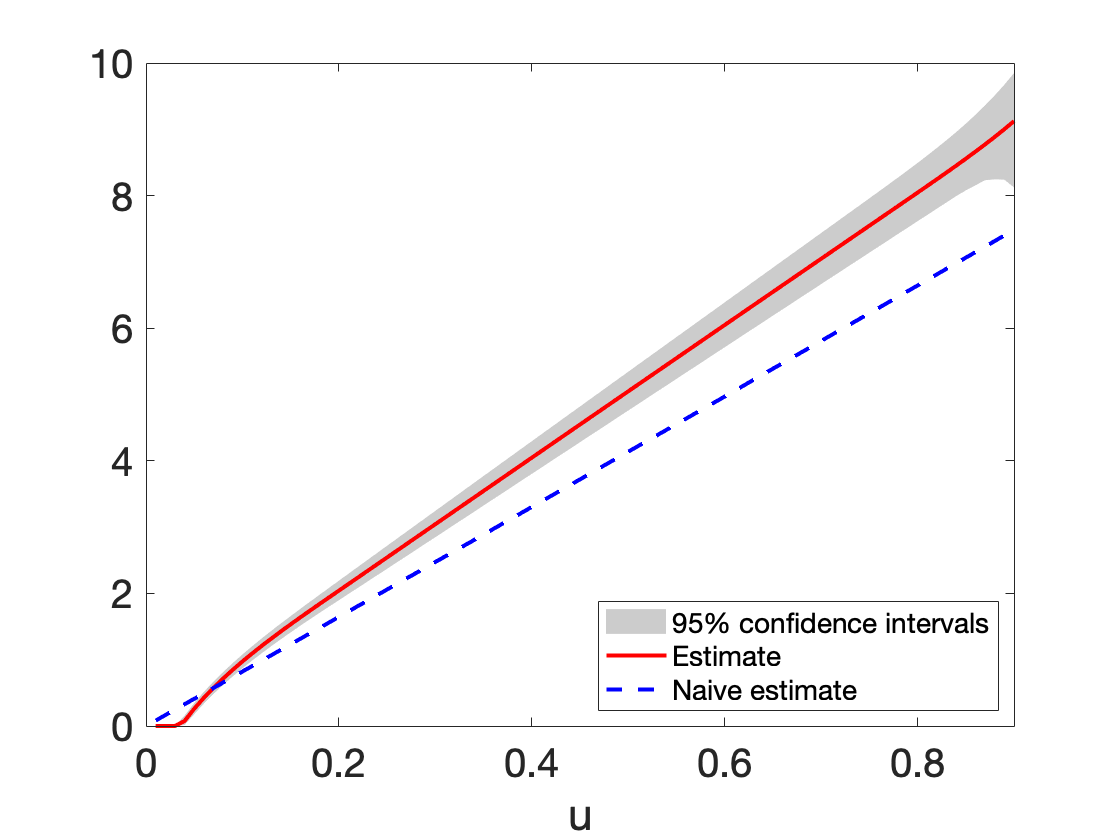}
  \caption{Average of the estimated regression function for $z=0$.}
   \label{fig:tracephi0} 
    \end{minipage}
    \quad\quad
    \begin{minipage}{0.48\textwidth}
        \centering
   \includegraphics[width=80mm]{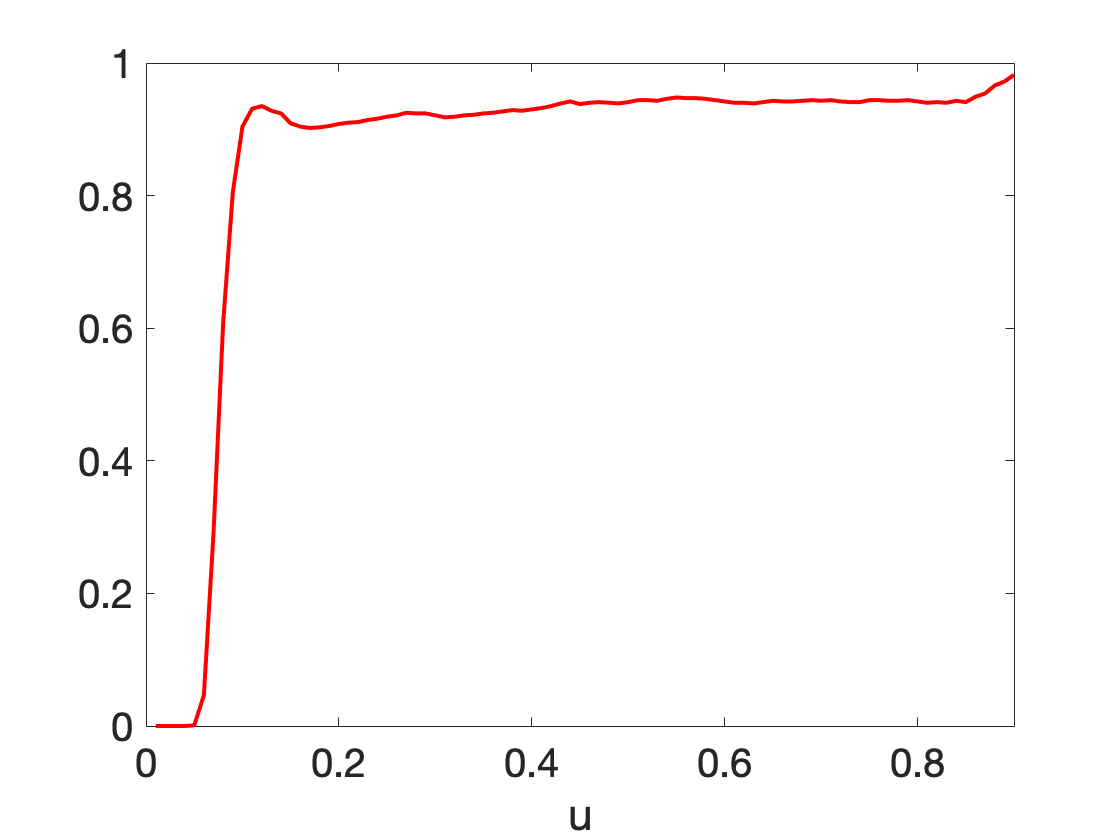}
  \caption{Coverage of confidence intervals for the regression function for $z=0$.}
      \label{fig:covphi0}
      \end{minipage}

\end{figure}

\begin{figure}[H]
\begin{minipage}{0.48\textwidth}
  \centering
  \includegraphics[width=80mm]{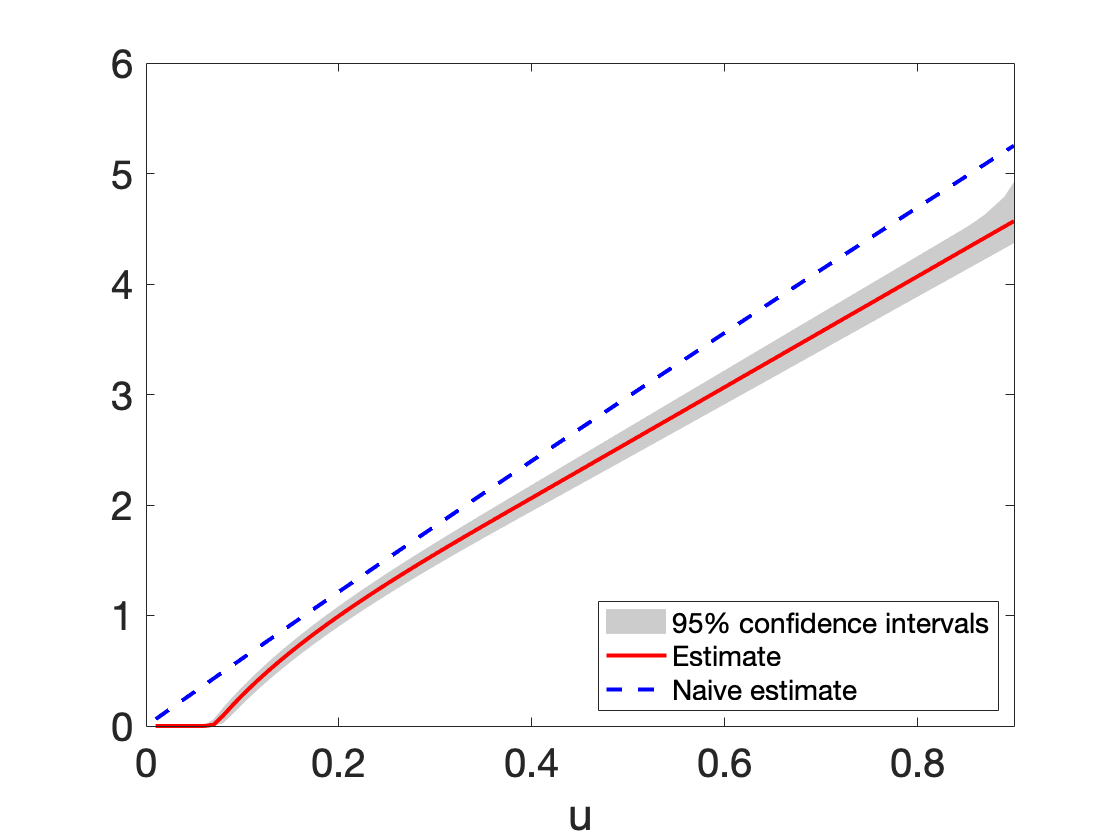}
  \caption{Average of the estimated regression function for $z=1$.}
   \label{fig:tracephi1} 
    \end{minipage}
    \quad\quad
    \begin{minipage}{0.48\textwidth}
        \centering
   \includegraphics[width=80mm]{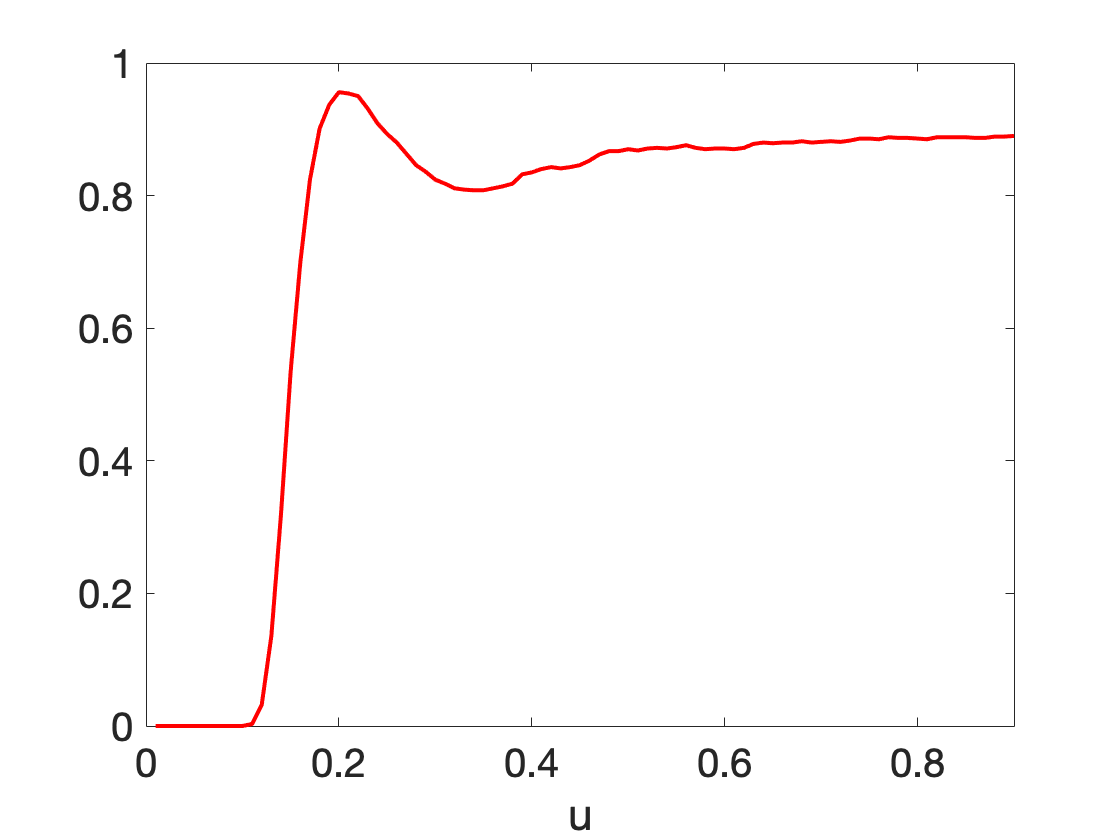}
  \caption{Coverage of confidence intervals for the regression function for $z=1$.}
      \label{fig:covphi1}
      \end{minipage}

\end{figure}

\begin{figure}[H]
\begin{minipage}{0.48\textwidth}
  \centering
  \includegraphics[width=80mm]{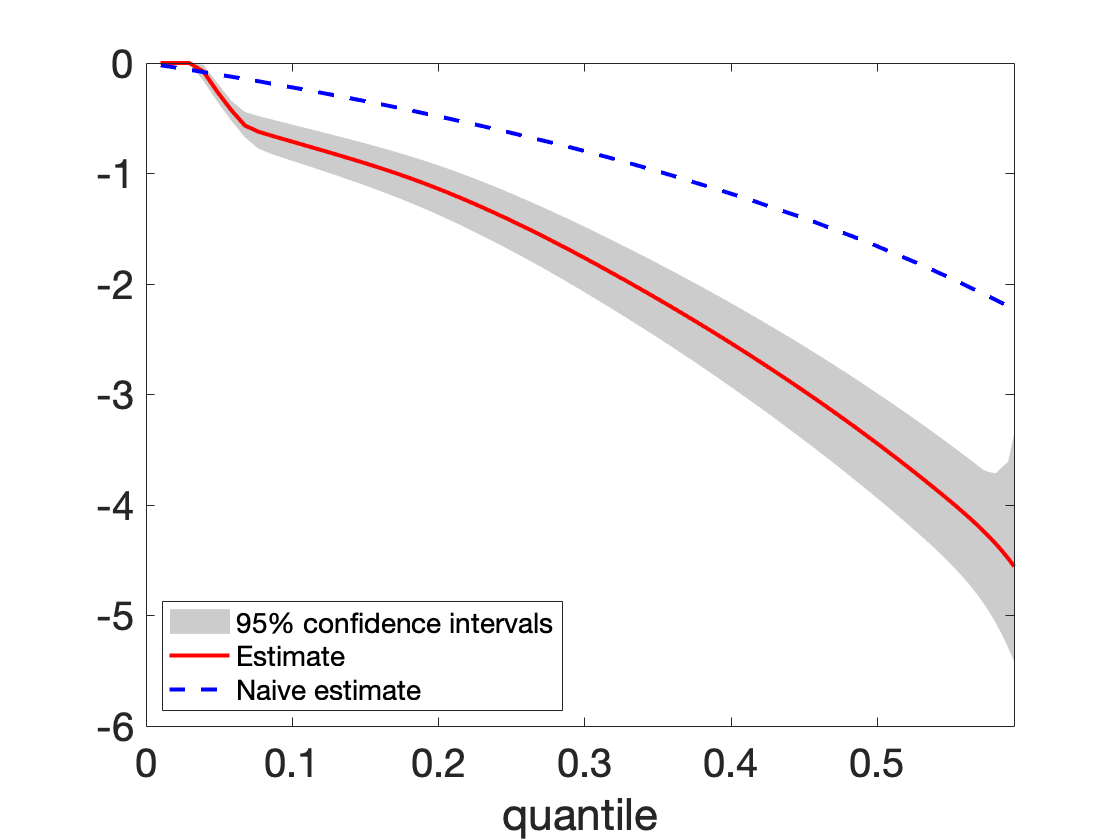}
  \caption{Average of the estimated QTE for different quantile levels.}
   \label{fig:traceQTE} 
    \end{minipage}
    \quad\quad
    \begin{minipage}{0.48\textwidth}
        \centering
   \includegraphics[width=80mm]{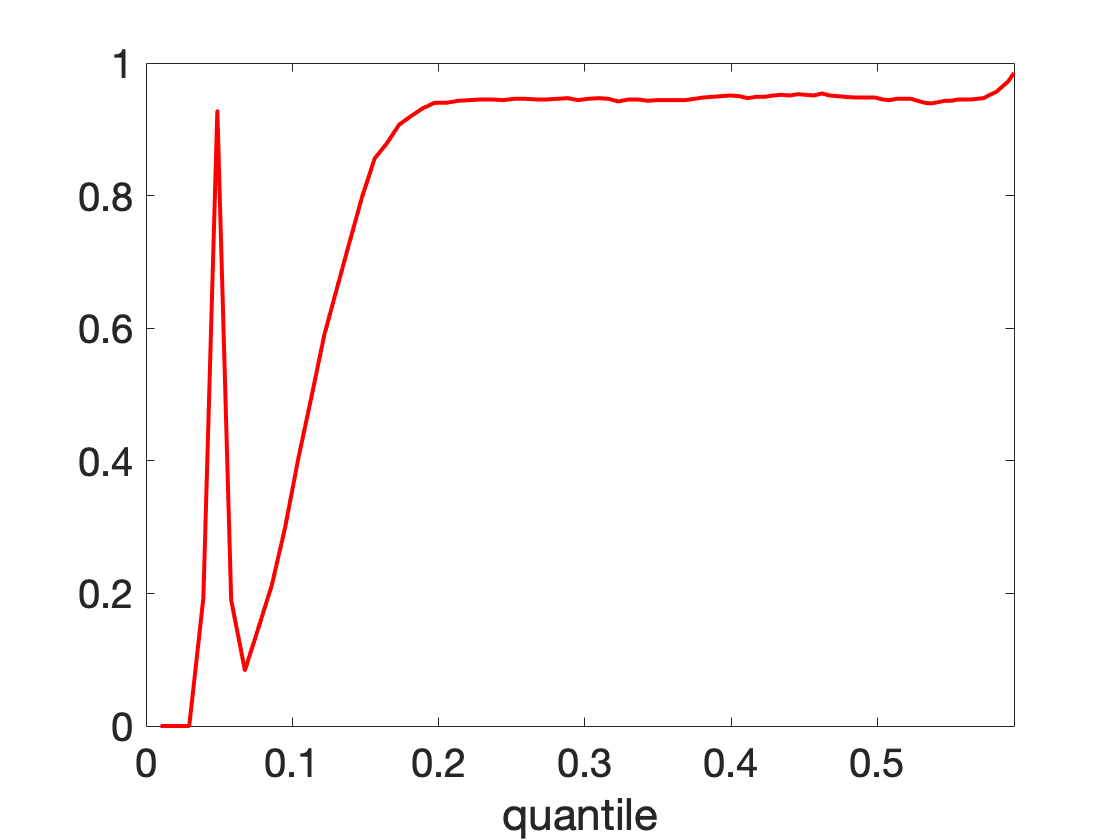}
  \caption{Coverage of confidence intervals for the QTE for different quantile levels.}
      \label{fig:covQTE}
      \end{minipage}

\end{figure}

 Figure \ref{fig:remainder_sim} reports the value of $\big\|A(\widehat \varphi, \widehat{S})(u) \big\|_{2}^2$ on these points. $\big\|A(\widehat \varphi, \widehat{S})(u) \big\|_{2}^2$ starts increasing slightly after $u$ reaches $0.9$. According to the rule of thumb of Section \ref{sec.pract}, we use the estimation results of Section \ref{sec.est} for $u$ below $0.9$. $\big\|A(\widehat \varphi, \widehat{S})(u) \big\|_{2}^2$ is large for $u$ close to zero because kernel estimators perform poorly near boundaries. Figure~\ref{fig:tracephi0} displays the estimated regression function $\widehat{\varphi}_0(\cdot)$. On the same graph, we also report the naive estimate of  $\varphi_0(\cdot)$ obtained by inversion of the Nelson-Aalen estimator of the cumulative hazard of $T$ given $Z=0$. The true regression function is omitted because it is undistinguishable from our estimate on the plot (except for very low values of $u$).
We also present the $95\%$ confidence intervals of  $\widehat{\varphi}_0(u)$ computed using $200$ bootstrap draws with the approach described in Section 4.4. Figure~\ref{fig:covphi0} contains the coverage of these confidence intervals. Figures ~\ref{fig:tracephi1} and ~\ref{fig:covphi1} report the same information for 
$\widehat{\varphi}_1(\cdot)$. The estimated quantile treatment effects and corresponding confidence intervals are presented in Figures ~\ref{fig:traceQTE} and ~\ref{fig:covQTE}. The quantile level $1-e^{-u}$ ranges from 0 to 0.6.

 It is clear that the naive estimator (in dashed blue) is biased. On the contrary, our estimator recovers precisely the shape of the true regression function. The confidence intervals have almost nominal coverage except for very low values of $u$ where $S$ is poorly estimated because of boundary properties of kernel estimators. In the supplementary material, we present the results of the simulations when a local polynomial of degree one is used to smooth the Kaplan-Meier estimator of the survival function. The coverages of the confidence intervals with local polynomial smoothing are much closer to $0.95$ for very low values of $u$. This is  because local polynomials do not suffer from the bias of kernel estimators near boundaries. 

\subsection{Empirical application}

Our empirical application revisits the so-called Illinois Reemployment Bonus Experiment.  We refer to \url{https://www.upjohn.org/data-tools/employment-research-data-center/} \url{illinois-unemployment-incentive-experiments} for more details. Taking place between mid-1984 and mid-1985, this randomized control trial aimed at evaluating the effects of bonuses to employers or job seekers on unemployment duration. Specifically, eligible unemployment insurance beneficiaries were randomly allocated to one of the following groups: 
\begin{itemize}
\item[(i)] \textbf{The Job Search Incentive Experiment group (JSIE).} Job seekers of this group were eligible to a cash bonus of 500\$ if they were to find a job of at least 30h/week within 11 weeks from the beginning of their unemployment spell and held that job for 4 months. In order to receive the bonus, job seekers had to read the experiment description and sign an agreement form at the beginning of their unemployment spell;
\item[(ii)]\textbf{The Hiring Incentive Experiment group (HIE).} Companies hiring job seekers within 11 weeks from the beginning of their unemployment spell for a job of more than 30h/week lasting at least 4 months were eligible to a cash bonus of 500\$. In order for their employer to receive the bonus, job seekers had to read the experiment description and sign an agreement form at the beginning of their unemployment spell;
\item[(iii)]\textbf{The control group.} \end{itemize}

The following specifics are worth noting. First, to be part of one of the three aforementioned groups, claimants needed to be between 20 and 55 years old, and have a valid unemployment insurance (UI) claim.  Second, the unemployment duration was recorded as the number of weeks during which participants received unemployment benefits, which implies that the data are discrete. In this particular case, the true unemployment duration is continuous but is only observed by interval. This setting creates another source of partial identification studied for instance in \citet{manski2002inference}. Because an extension to such a setting is outside of the scope of this paper, we choose to neglect this issue and treat the data as continuous. Third, given that UI is granted for 26 weeks, we can only observe individuals up to the end of their UI claim, that is the data are right censored at 26 weeks. Around 21\% of the observations are censored.

To clarify the relevance of our model to the experiment, let us assume that we want to evaluate the causal effect of the cash bonus of the JSIE (or HIE) experiment. Although the allocation to the groups is randomly chosen, the decision to agree to participate in the experiment is endogenous. Let us consider a framework with two treatments as described by the following variable $Z$:
$$Z=\left\{\begin{array}{cl} 
2 &\text{if the individual is in the HIE group and agrees to participate}\\
1 &\text{if the individual is in the JSIE group and agrees to participate}\\
0 &\text{if the individual is assigned to the control group or refuses to participate.}\end{array} \right.$$
In order to deal with the selection bias, our strategy consists of using the group assignment as an instrument. Let us define
$$W=\left\{\begin{array}{cl} 
2 &\text{if the individual is assigned to the HIE group}\\
1 &\text{if the individual is assigned to the JSIE group}\\
0 &\text{if the individual is assigned to the control group.}\end{array} \right.$$\\
We report in Table \ref{tab2} the sample sizes for each combination of values of the couple $(W,Z)$. The refusal rates are respectively 16\% and 35\% in the JSIE and HIE experiment, which suggests large selection biases. In total, there are $12,101$ observations.

 \begin{table}[!ht]
    \center
    \begin{tabular}{|c|c|c|c|}
    \hline 
Sample size & $W=0$ &  $W=1$ &$W=2$\\
\hline
$Z=0$ & $3952$ &$1377$ &$659$\\
\hline 
$Z=1$ & $0$ & $2586$&$0$\\
\hline 
$Z=2$ & $0$ &$0$&$3527$\\
\hline
    \end{tabular}
    \caption{Sample sizes in the experiment.}
    \label{tab2}
\end{table}

This dataset has been extensively studied in the literature. Comparing `naively' treatment groups to the control group, \citet{WS} concluded that the JSIE bonus reduces (statistically significantly) the unemployment duration but found no significant evidence of an effect of the HIE bonus. Using a proportional hazards model with $W$ as explanatory variable, \citet{M88} reached the conclusion that the JSIE bonus increases the hazard rate. In \citet{BR}, the hazard rate follows a mixed proportional hazards model which is semiparametric. The authors develop a two-stage instrumental variable estimator, and find a stronger positive effect of the JSIE experiment than \citet{M88}. 

\subsubsection{Estimation under exact identification}

\begin{figure}[ht]
\centering
\includegraphics[width=80mm]{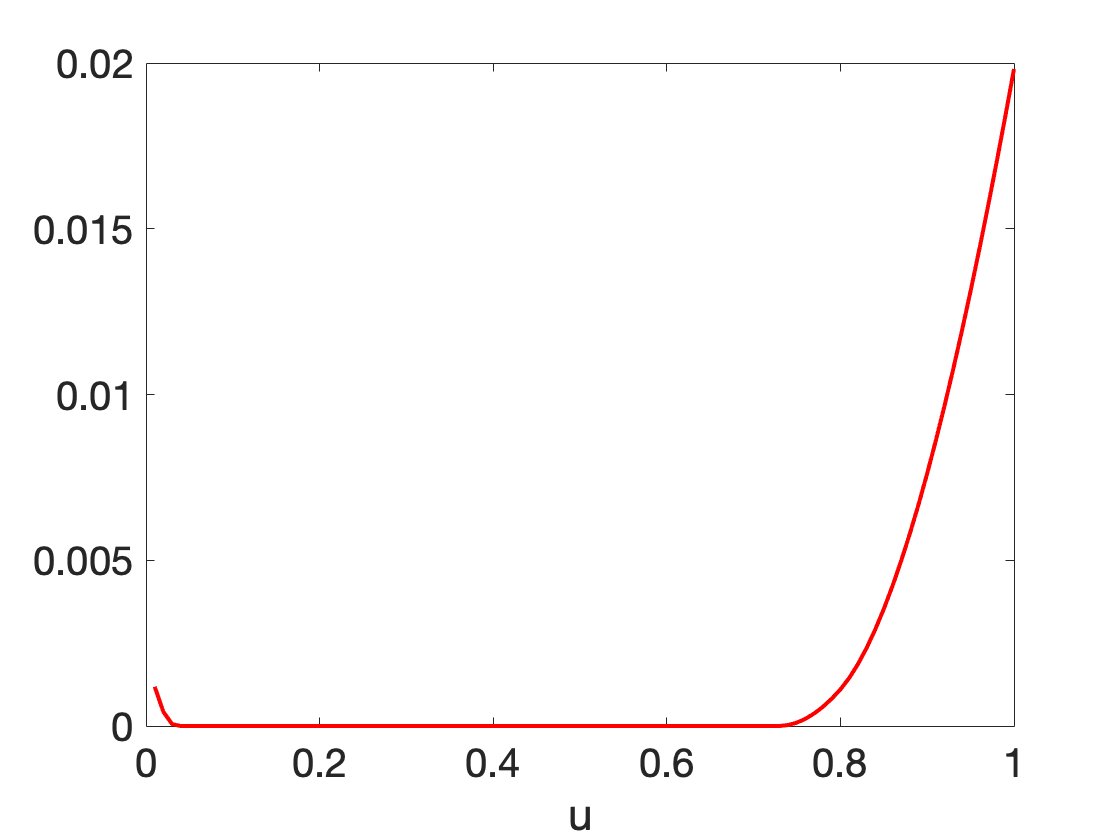}
\includegraphics[width=80mm]{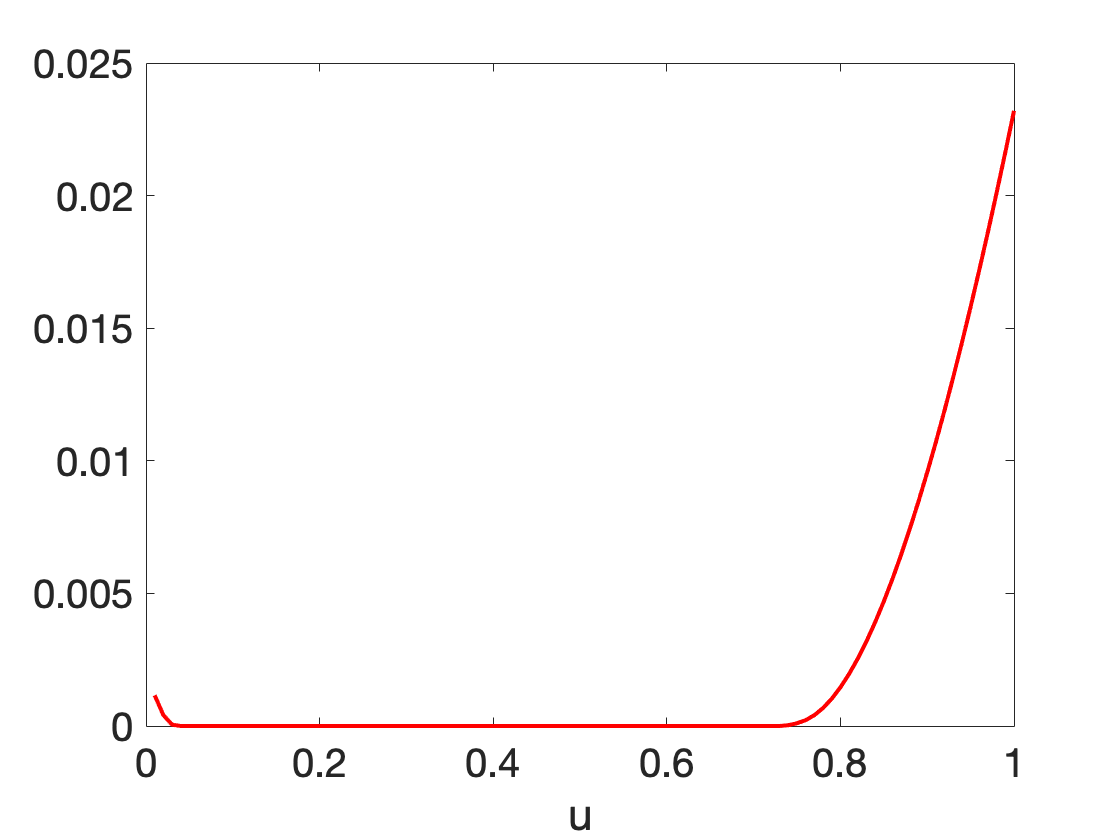}
\caption{Plot of $u\mapsto \big\|A(\widehat \varphi, \widehat{S})(u) \big\|_{2}^2$ when the analysis is restricted to the population for which $Z=0$ or $Z=1$ on the left and $Z=0$ or $Z=2$ on the right.}
\label{fig:remainder_emp}
\end{figure}

We estimate quantile treatment effects of the two different bonuses using the method described in this paper. An Epanechnikov kernel is used to smooth the survival function. The bandwidth is equal to $2$, although the results are not sensitive to the choice of the bandwidth as long as it is sufficiently larger than $1$. $\bar{T}$ is equal to $26$, that is the censoring duration. We compute the estimators of $\varphi$ on a grid for $u$ between $0.01$ and $1$ with step size $0.01$. Figure \ref{fig:remainder_emp} reports the value of $\big\|A(\widehat \varphi, \widehat{S})(u) \big\|_{2}^2$ on the grid, when the analysis is restricted to the population for which $Z=0$ or $Z=1$ on the left and $Z=0$ or $Z=2$ on the right. It can be seen that $\big\|A(\widehat \varphi, \widehat{S})(u) \big\|_{2}^2$ starts increasing slightly after the value of $u$ reaches $0.7$. Applying the analysis of Section \ref{sec.pract}, we choose to use the estimation results of Section \ref{sec.est} for $u$ below $0.7$. The corresponding quantiles $1-e^{-u}$ are between 0 and around 0.5. Note that $\big\|A(\widehat \varphi, \widehat{S})(u) \big\|_{2}^2$ is large when $u$ is close to zero, this is because kernel estimators perform poorly near boundaries.

\begin{figure}[ht]
\centering
\includegraphics[width=80mm]{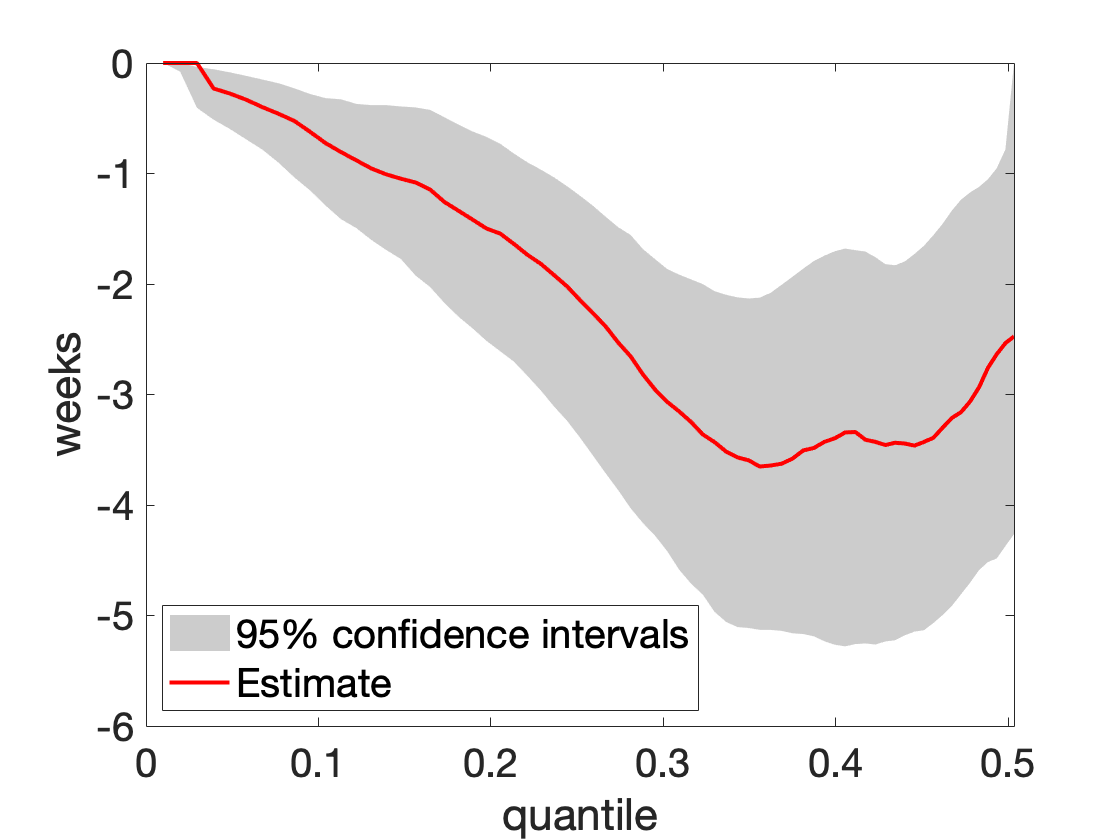}
\includegraphics[width=80mm]{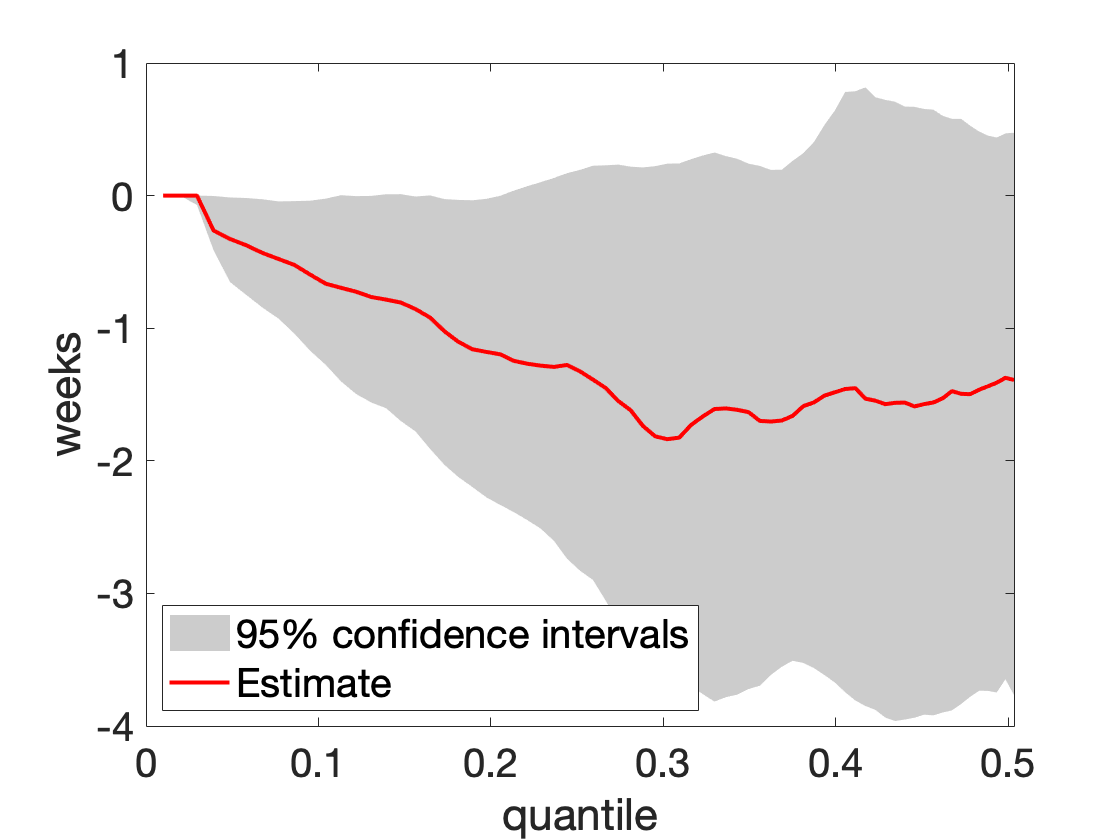}
\caption{Estimated quantile treatment effects and confidence intervals for $Z=0$ and $Z=1$ (left) and for $Z=0$ and $Z=2$ (right) using kernel smoothing.}
\label{fig:comp3}
\end{figure}

Figure \ref{fig:comp3} reports estimated quantile treatment effects on this range. We estimate $95\%$ confidence intervals using $1000$ bootstrap draws. Consistent with previous studies, we conclude that the JSIE bonus has a stronger effect than the HIE bonus and both reduce the unemployment duration. For most values of $u$ in the selected range, the quantile treatment effect of the JSIE bonus is significant while the one of the HIE bonus is not. Estimation results using a local polynomial of degree one are similar and, hence, are not reported.

\subsubsection{Partial identification}
Let us now apply the partial identification results from Section \ref{sec.id}. This empirical application fits in the case of the triangular system of equations of Section \ref{sec.speid} if we restrict the data to individuals for which $W=0$ or $W=1$. Then, we can estimate the outer set of the identified set for values of $u$ for which the following system does not have an exact solution $(\theta_1,\theta_2)\in[0,c_0)\times [0,c_0)$, where $c_0=26$:
\begin{eqnarray} \label{trian2}
\left\{
\begin{array}{l}
\widehat{S}(\theta_1 ,0|0) = e^{-u} \\
\widehat{S}(\theta_2,1|1) + \widehat{S}(\theta_1,0|1) = e^{-u}.
\end{array}
\right. 
\end{eqnarray}

Let us define $u(c_0)=-\log(\widehat{S}(c_0 ,0|0) )$, which is equal to around 0.73. For $u\ge u(c_0)$ it is clear that the system \eqref{trian2} has no solution in $[0,c_0)\times [0,c_0)$. We are, hence, in the second case described in Section \ref{sec.speid}, and we can estimate an outer set to the identified set for $(\varphi_0(u(c_0)), (\varphi_1(u(c_0)))^\top$  by $[c_0,\infty)\times [0, \widehat{\theta}_2]$, where $\widehat{\theta}_2$ is defined by $\widehat{S}(\widehat{\theta}_2,1|1)= e^{-u}- \widehat{S}(c_0,0|1)$. For $u=u(c_0)$, we obtain $\widehat{\theta}_2=23.31$ and the estimated identified set is $[26,\infty)\times [0, 23.31]$. This implies that the estimated  outer set of the quantile treatment effect of the JSIE bonus at $u(c_0)$ is $(-\infty,-2.69]$, which suggests that the true quantile treatment effect is lower than $-2.69$.


\section{Conclusion and future research}
\label{sec:conc}

In this paper, we studied the issue of identification and estimation of an endogenous categorical regressor on a duration, in the presence of a categorical instrumental variable.  We developed partial, local and global identification results in the presence of random right censoring in a fully nonparametric framework. The causal effect is expressed as the solution of a nonlinear system of equations. We derived the asymptotic properties of the solution to an estimated version of this system. Our simulations exhibit excellent performance of our approach in finite samples. We showed how to revisit an empirical application using our methodology.

Several extensions of our model may be interesting. One could consider the case where $Z$ and/or $W$ are continuous variables, or even dynamic processes. For instance, $Z$ could represent a treatment that can be obtained throughout the unemployment spell of an individual, e.g. a job training. \red{Specifically, if the variables $Z$ and $W$ are continuous (but not time dependent), equation \eqref{Sy2} is replaced by a nonlinear integral equation \begin{equation}\label{continuouS}\int S(\varphi(z,u),z\vert w)dz=e^{-u},\end{equation}
where $\varphi$ is the functional parameter of interest and the function $S$ may be nonparametrically estimable. This equation generates an ill-posed problem and its
resolution requires regularization (see e.g. \citet{kaltenbacher2008iterative} or for econometric applications, \citet{C}). The most convenient
solution is to implement a sequential resolution algorithm stopped in a suitable way. The
case with censoring and partial identification raises an original question in this
framework. The extension of our approach to dynamic treatments (and possibly dynamic 
instruments) would be based on a generalisation of Lemma \ref{hazard}. For example, $Z$ may be replaced by $Z_t = I(t \ge \epsilon)$, where $\epsilon$ is the random starting time of the treatment. Dynamic
treatments also generate new censoring mechanisms and will be considered in future
research. These problems are discussed in the unpublished working paper by \citet{florens2010endogeneity}. }

Another valuable generalization would allow for additional covariates $X$. If the latter are discrete, they may be included in the variable $Z$ or the econometrician could use the approach of the present paper conditional on $X$. If they are continuous, the problem is more complicated. If $T=\varphi(Z,X,U)$, analogously to \eqref{continuouS}, we would have the system 
$$ \int S(\varphi(z,x,u),z,x\vert w)dz=e^{-u}$$ or the equivalent with a sum if $Z$ is discrete. One would need to regularize the estimator as discussed in the previous paragraph.

As shown by the empirical application, an interesting extension would concern discrete outcomes. If a discrete outcome arises because $T$ is continuous but is only observed as an interval, then we have a second source of partial identification. Such a censoring mechanism has been studied in the partial identification literature, see for instance \citet{manski2002inference}. If instead, the true outcome is discrete, then we need to allow $\varphi_Z(\cdot)$ to be a step function. In this case \eqref{Sy} cannot be obtained because $\varphi_Z(\cdot)$ would not be strictly increasing.

Another possible source of partial identification is nonrandom censoring. Such a case was studied in \citet{blanco2019bounds} using principal stratification. A similar approach may be tried in the context of this paper.

Finally, inference in the partially identified case due to censoring is an appealing topic.  We conjecture that a bootstrap approach may be used. For each bootstrap simulation, the function $\min_k R_{k,u}(\theta)$ may be evaluated and then the corresponding outer set can be computed. Then, for any point outside of $[0,c_0)^L$, one may compute the bootstrap probability to belong to this outer set and confidence regions of the outer set may be derived. Another possibility would be to use results on models based on conditional moment inequalities such as in \citet{andrews2013inference}.

\bibliographystyle{dcu}
\bibliography{NPIV_CENSORING}

\appendix

\renewcommand{\theequation} {A.\arabic{equation}}

\section{Appendix A: Proofs of the results of Sections 2, 3 and 4}
\label{sec.A}

\noindent
{\bf Proof of Lemma \ref{hazard}.}
Since $U \sim \mbox{Exp}(1)$ and $U$ and $W$ are independent, we have for $w=w_1,\dots, w_K$ and $u\in \mathbb{R}_+$ that
$$ \sum_{\ell=1}^L f_{U,Z}(u,z_\ell|w) = e^{-u} \quad \mbox{ and } \quad \sum_{\ell=1}^L S_{U,Z}(u,z_\ell|w) = e^{-u}, $$
where $S_{U,Z}(u,z|w) = \mathbb{P}(U \ge u,Z=z|W=w)$ and $f_{U,Z}(u,z|w) = -\frac{\partial S_{U,Z}}{\partial u}(u,z|w)$.  Hence,
$$ \sum_{\ell=1}^L  \varphi_{z_\ell}'(u) f(\varphi_{z_\ell}(u),z_\ell|w) = e^{-u} \quad \mbox{ and } \quad \sum_{\ell=1}^L S(\varphi_{z_\ell}(u),z_\ell|w) = e^{-u}, $$
where $S(t,z|w) = \mathbb{P}(T \ge t,Z=z|W=w)$ and $f(t,z|w) = -\frac{\partial S}{\partial t}(t,z|w)$.  This implies that
\begin{eqnarray*}
1 &=& \sum_{\ell=1}^L \varphi_{z_\ell}'(u)  \frac{f(\varphi_{z_\ell}(u),z_\ell|w)}{\sum_{\ell'=1}^L S(\varphi_{z_{\ell'}}(u),z_{\ell'}|w)} \\
&=& \sum_{\ell=1}^L \varphi_{z_\ell}'(u)  \frac{f(\varphi_{z_\ell}(u),z_\ell|w)}{S(\varphi_{z_\ell}(u),z_\ell|w)} \frac{S(\varphi_{z_\ell}(u),z_\ell|w)}{\sum_{\ell'=1}^L S(\varphi_{z_{\ell'}}(u),z_{\ell'}|w)} \\
&:=& \sum_{\ell=1}^L \varphi_{z_\ell}'(u) R_1(u,z_\ell|w) R_2(u,z_\ell|w),
\end{eqnarray*}
say.  Note that, for $z=z_1,\dots,z_L$ and $w=w_1,\dots,w_K$, we have
$$ R_1(u,z|w) = \frac{f(\varphi(z,u),z|w) \mathbb{P}(Z=z|W=w)}{S(\varphi(z,u),z|w) \mathbb{P}(Z=z|W=w)} = \frac{f(\varphi(z,u)|z,w)}{S(\varphi(z,u)|z,w)} = h(\varphi(z,u)|z,w) $$
and $R_2(u,z|w) = p(z|U \ge u, w)$, since 
$$ p(z|U \ge u, w) = \frac{\mathbb{P}(Z=z|W=w) S_U(u|z,w)}{\sum_{\ell=1}^L \mathbb{P}(Z=z_\ell|W=w) S_U(u|z_\ell,w)} = \frac{S(\varphi(z,u),z|w)}{\sum_{\ell=1}^L S(\varphi_{z_\ell}(u),z_\ell|w)}. $$
Hence, integrating over $u$ (and replacing $u$ by $v$) yields
\begin{eqnarray*}
u &=& \int_0^u \sum_{\ell=1}^L\varphi_{z_\ell}'(v)  h(\varphi(z_\ell,v)|z_\ell,w) p(z_\ell|U \ge v, w) dv \\
&=& \sum_{\ell=1}^L \int_0^{\varphi_{z_\ell}(u)} h(s|z_\ell,w) p(z_\ell|T \ge s, w) \, ds.
\end{eqnarray*}
\vspace*{-1.5cm}

\hfill $\Box$\\

\bigskip 

\noindent
{\bf Proof of Theorem \ref{C}.} We have 
\begin{align} 
\notag 
\|A(\widehat{\varphi},S)\|&\le \|A(\widehat{\varphi},\widehat{S})-A(\widehat{\varphi},S)\|+ \|A(\widehat{\varphi},\widehat{S})\|\\
\label{ineqM}&= O_P(r_n)+\|A(\widehat{\varphi},\widehat{S})\|.
\end{align} 
Now, notice that, by definition,
\begin{align*}
\|A(\widehat{\varphi},\widehat{S})\|
&\le \frac{1}{\sigma_{\min}(V)}\|A(\widehat{\varphi},\widehat{S})\|_V \\
&= \frac{1}{\sigma_{\min}(V)}\inf_{\theta \in\mathcal{F}_Z^{\bar{U},\bar{T}}}\|A(\theta,\widehat{S})\|_V\\
&\le  \frac{1}{\sigma_{\min}(V)}\left( \inf_{\theta \in \mathcal{F}_Z^{\bar{U},\bar{T}}}\|A(\theta,S)\|_V+ \sup_{\theta \in \mathcal{F}_Z^{\bar{U},\bar{T}}}\|A(\theta,\widehat{S})-A(\theta,S)\|_V\right)\\
&\le  \frac{1}{\sigma_{\min}(V)}\left( \inf_{\theta \in \mathcal{F}_Z^{\bar{U},\bar{T}}}\|A(\theta,S)\|_V+ \sigma_{\max}(V)\sup_{\theta \in \mathcal{F}_Z^{\bar{U},\bar{T}}}\|A(\theta,\widehat{S})-A(\theta,S)\|\right)\\
&= 0+O_P(r_n)\quad \text{(By assumptions (V) and (C))},
\end{align*}
where $\sigma_{\min}(V)$ (resp. $\sigma_{\max}(V)$) is the lower bound (resp. upper bound) defined in assumption (V).

This and \eqref{ineqM} imply that $\|  A(\widehat{\varphi},S)\|=O_P(r_n)$, which yields that $\| \widehat{\varphi}-\varphi \|=o_P(1)$ using Assumption (C)(i). Hence, by (C)(ii), $\| \widehat{\varphi}-\varphi \| \le c\|  A(\widehat{\varphi},S)\|$ with probability approaching $1$, which shows the result. \hfill $\Box$ \\[-.3cm]

\noindent
{\bf Proof of Theorem \ref{AS}.} The first-order conditions of the optimisation program \eqref{prog} are 
$$\Gamma(\widehat{\varphi},\widehat{S})^\top V A(\widehat{\varphi},\widehat{S})=0,$$
which leads to 
\begin{align*}
& \Gamma(\widehat{\varphi},\widehat{S})^\top VA(\varphi,\widehat{S}) 
+ \Gamma(\widehat{\varphi},\widehat{S})^\top V\Gamma(\varphi,\widehat{S})(\widehat{\varphi}-\varphi)
+ \Gamma(\widehat{\varphi},\widehat{S})^\top V R=0,
\end{align*}
where $R= A(\widehat{\varphi},\widehat{S})-A(\varphi,\widehat{S})-\Gamma(\varphi,\widehat{S})(\widehat{\varphi}-\varphi)$.
Then, 
\begin{align*}
\sqrt{n} (\widehat{\varphi}-\varphi)
= & -[ \Gamma(\widehat{\varphi},\widehat{S})^\top V\Gamma(\varphi,\widehat{S})]^{-1}\Gamma(\widehat{\varphi},\widehat{S})^\top V\sqrt{n}A(\varphi,\widehat{S})\\
& - [ \Gamma(\widehat{\varphi},\widehat{S})^\top V\Gamma(\varphi,\widehat{S})]^{-1} \Gamma(\widehat{\varphi},\widehat{S})^\top V\sqrt{n} R.
\end{align*}
Now, by Assumption (N) and Theorem \ref{C}, $\sqrt{n}R=O_P(\sqrt{n}r_n^2)=o_P(1)$, and hence 
$$[ \Gamma(\widehat{\varphi},\widehat{S})^\top V\Gamma(\varphi,\widehat{S})]^{-1} \Gamma(\widehat{\varphi},\widehat{S})^\top V \sqrt{n} R=o_P(1).$$ We conclude using the continuous mapping theorem and Assumption (N). \hfill $\Box$

\section{Appendix B: Kernel smoothing of the Kaplan-Meier estimator}
 
In this section, we discuss sufficient conditions for Assumptions (C) and (N). Let $\widehat{S}$ be defined as in \eqref{KM}. We introduce the following hypothesis and the subsequent lemma.
\begin{itemize}
\item[\textbf{(K)}] 
\begin{itemize} 
\item[(i)] $T$ has a continuous distribution conditional on the event $\{Z=z,W=w\}$ for all $z$ and $w$. Its conditional density has a derivative which is bounded on $[0,\bar{T}]$; 
\item[(ii)] $K$ is a bounded and differentiable probability density function, it has mean zero and a bounded support; 
\item[(iii)] $n\epsilon^4 \to 0$;
\item[(iv)] The Fr\'echet differential $\Gamma(\varphi,S)$ is invertible;
\item[(v)] There exists $\xi>0$ for which $S(\bar T,z|w)/\mathbb{P}(Z=z|W=w)>\xi$ for all $z\in\{z_1,\dots, z_L\}$ and $w\in\{w_1,\dots, w_K\}$ such that $\mathbb{P}(Z=z|W=w)>0$.
\end{itemize}
\end{itemize}

\begin{Lemma} \label{K}
Under Assumptions (V), (C)(i), (C)(ii), (K) and if the mapping $\Gamma(\varphi,S)^\top V\Gamma(\varphi,S)$ is invertible on $[0,\bar{U}]$, then $\widehat{S}$  satisfies Assumptions (C)(iii) with $r_n=n^{-1/2}$ and (N) and the bootstrap approximation of Section 4.4 works in the sense that $\sqrt{n}(\widehat{\varphi}_b-\widehat{\varphi})$ converges weakly to the same mean zero Gaussian process as $\sqrt{n}(\widehat{\varphi}-\varphi)$. 
\end{Lemma}

The lemmas below constitute the proof of Lemma \ref{K}.
We introduce the following notations. 
For $h\in \mathcal{F}_{Z,W}$, we define 
$$ \vert\vert h \vert\vert_{\infty}=\underset{t\in [0,\bar{T}], z=z_1,\dots,z_L, w=w_1,\dots,w_K} {\sup}\vert h(t, z,w) \vert. $$ 
Also, for $t \in [0,\bar{T}], z=z_1,\ldots,z_L, w=w_1,\ldots,w_K$, let $\widetilde{S}(t|z,w)=\int S(t-s\epsilon|z,w)K(s)ds$.

\begin{Lemma}\label{FAN}
The processes $\sqrt{n}\{\widehat{S}(t,z|w)-p_{zw}\widetilde{S}(t|z,w)\}$ and $\sqrt{n}(\widehat{S}_b(t,z|w)-\widehat{S}(t,z|w))$, $t \in [0,\bar{T}], z=z_1,\ldots,z_L, w=w_1,\ldots,w_K$, converge weakly to the same mean zero Gaussian process.
\end{Lemma}
{\bf Proof.}  By the Donsker theorem, the process 
$$ \frac{1}{\sqrt{n}} \left\{\left(\begin{array}{c}N_{z,w}(t)\\Y_{z,w}(t)\\ Y_{z,w}\\  Y_w\end{array}\right)-\left(\begin{array}{c}\mathbb{P}(Y\le t,Z=z,W=w,\delta=1)\\ \mathbb{P}(Y\ge t,Z=z,W=w)\\ \mathbb{P}(Z=z,W=w)\\\mathbb{P}(W=w) \end{array}\right)\right\}, $$
$t\in[0,\bar{T}], z=z_1,\ldots,z_L, w=w_1,\ldots,w_K$, converges weakly to a mean zero Gaussian process. Note that the mapping from this process to $\widehat{S}$ is Hadamard differentiable as it is a composition of Hadamard differentiable mappings. Indeed, the quotient is Hadamard differentiable and condition K(v) guarantees that the mapping yielding the Kaplan-Meier estimator in \eqref{KM} by Lemma 20.14 in \citet{VdV}. As a linear operator, the mapping embodying kernel smoothing as in \eqref{smooth} is Hadamard differentiable. Finally, the product is Hadamard differentiable. We obtain that the process $\sqrt{n}\{\widehat{S}(t,z|w)-p_{zw}\widetilde{S}(t|z,w)\}$, $t \in [0,\bar{T}], z=z_1,\ldots,z_L, w=w_1,\ldots,w_K$, converges to a mean zero Gaussian process using the functional Delta-method. 
Similarly, using the bootstrap functional Delta-method as in Theorem 3.9.11 in \citet{VdVW}, we find that $\sqrt{n}(\widehat{S}_b(t,z|w)-\widehat{S}(t,z|w))$, $t \in [0,\bar{T}], z=z_1,\ldots,z_L, w=w_1,\ldots,w_K$, converges weakly to the same mean zero Gaussian process. \hfill $\Box$\\[-.3cm]

\begin{Lemma} \label{true} 
Under Assumption (K), it holds that
$$\underset{t\in [0,\bar{T}], z=z_1,\dots,z_L, w=w_1,\dots,w_K} {\sup} |p_{zw}\widetilde{S}(t|z,w)-S(t, z| w)|
=o\Big(\frac{1}{\sqrt{n}}\Big).$$
\end{Lemma}
{\bf Proof.}  From a Taylor expansion, as the second derivative of $S(\cdot|z,w)$ is bounded on $[0,\bar{T}]$, we have 
$$ S(t-s \epsilon \vert z,w) = S(t\vert z,w ) + s\epsilon\frac{\partial S}{\partial t}(t\vert z,w) +O(\epsilon^2), $$
uniformly in $t\in [0,\bar{T}], z=z_1,\dots,z_L, w=w_1,\dots,w_K$ and $s$ in the support of $K$.  Therefore, we obtain
\begin{align*}
& \sup_{t\in [0,\bar{T}], z=z_1,\dots,z_L, w=w_1,\dots,w_K}\Big|\int \Big(S(t-s \epsilon \vert z,w)-S(t\vert z,w )-s\epsilon\frac{\partial S}{\partial t}(t\vert z,w)\Big)K(s)ds\Big|\\
&=\sup_{t\in [0,\bar{T}], z=z_1,\dots,z_L, w=w_1,\dots,w_K}\Big|\int (S(t-s \epsilon \vert z,w)-S(t\vert z,w ))K(s)ds\Big|\\
& =  O(\epsilon^2),
\end{align*}
where the last sequence of equalities is a consequence of the fact that $K$ has mean zero, is bounded and has bounded support. This leads to 
$$ \sup_{t \in [0,\bar{T}], z=z_1,\ldots,z_L, w=w_1,\ldots,w_K} |\widetilde{S}(t|z,w)-S(t|z, w)|=O(\epsilon^2),$$
because $\frac{\partial^2 S}{\partial t^2}(t\vert z,w)$ is bounded on $[0,\bar{T}]$, by Assumption (K). This yields that 
\begin{align*}
& \sup_{t \in [0,\bar{T}], z=z_1,\ldots,z_L, w=w_1,\ldots,w_K} |p_{zw}\widetilde{S}(t|z,w)-S(t, z| w)| \\
 &\le \sup_{t \in [0,\bar{T}], z=z_1,\ldots,z_L, w=w_1,\ldots,w_K} |\widetilde{S}(t|z,w)-S(t|z, w)|\\
 &=O(\epsilon^2)=o \Big(\frac{1}{\sqrt{n}}\Big)\quad \text{(by Assumption (K))}.
\end{align*}
\vspace*{-1.5cm} 

\hfill $\Box$\\[-.3cm]
 
\begin{Lemma} \label{finalcv} 
Under Assumption (K), it holds that $\| \widehat{S} -S\|_{\infty}=O_P(\frac{1}{\sqrt{n}})$ and $\| \widehat{S}_b -\widehat{S}\|_{\infty}=O_P(\frac{1}{\sqrt{n}})$.
\end{Lemma}
{\bf Proof.}  By Lemma \ref{FAN} and the continuous mapping theorem, we have 
$$\underset{t\in [0,\bar{T}], z=z_1,\dots,z_L, w=w_1,\dots,w_K} {\sup}\vert \widehat{S}(t,z|w)-p_{zw}\widetilde{S}(t|z,w) \vert=O_P \Big(\frac{1}{\sqrt{n}}\Big).$$
As $ | \widehat{S}(t,z|w) -S(t, z\vert w)|\le \vert \widehat{S}(t,z|w)-p_{zw}\widetilde{S}(t|z,w) \vert+|p_{zw}\widetilde{S}(t|z,w)-S(t, z| w)|$, the first result is a consequence of Lemma \ref{true}. 
The fact that $\| \widehat{S}_b -\widehat{S}\|_{\infty}=O_P(\frac{1}{\sqrt{n}})$ is a direct corollary of Lemma \ref{FAN}.
 \hfill $\Box$\\
 
 \begin{Lemma}\label{cvS} Under Assumption (K), the process
$\sqrt{n}(\widehat{S}(t,z|w)-S(t,z|w))$, $t \in [0,\bar{T}], z=z_1,\ldots,z_L, w=w_1,\ldots,w_K$, converges weakly to a mean zero Gaussian process.
\end{Lemma}
{\bf Proof.}  Using Slutsky's theorem, Lemmas \ref{FAN} and \ref{true} and the fact that, for $t\in [0,\bar{T}], z=z_1,\dots,z_L, w=w_1,\dots,w_K$,
$$\widehat{S}(t,z|w) -S(t, z\vert w) =  \widehat{S}(t,z|w)-p_{zw}\widetilde{S}(t|z,w) +p_{zw}\widetilde{S}(t|z,w)-S(t, z| w),$$
we obtain that $\sqrt{n}(\widehat{S}-S )$ converges weakly to a mean zero Gaussian process.  \hfill $\Box$\\[-.3cm]
 
\begin{Lemma}\label{cvM} 
Under Assumption (K), the process
$\sqrt{n}A(\varphi,\widehat{S})(u)$, $u \in [0,\bar{U}]$, converges weakly to a mean zero Gaussian process.
\end{Lemma}
{\bf Proof.}  By Lemma \ref{cvS}, $\sqrt{n}(\widehat{S}-S )$ converges weakly to a mean zero Gaussian process. The mapping $h \mapsto A(\varphi,h)$, $h \in \mathcal{F}_{Z,W}$, is Hadamard differentiable (because it is linear), therefore, by the functional Delta-method, we have that $\sqrt{n}(A(\varphi,\widehat{S})-A(\varphi,S))=\sqrt{n}A(\varphi,\widehat{S})$ converges weakly to a mean zero Gaussian process. \hfill $\Box$\\[-.3cm]

\begin{Lemma}\label{cvderiv} 
Under Assumption (K),
$\Gamma(\varphi,\widehat{S})(u)$, $u \in [0,\bar{U}]$ exists and converges in probability to $\Gamma(\varphi,S)(u)$, $u \in [0,\bar{U}]$. Also, if $\widehat{\varphi}\xrightarrow{P}\varphi$, then $\Gamma(\widehat{\varphi},\widehat{S}_b)(u)$, $u \in [0,\bar{U}]$ exists and converges in probability to $\Gamma(\varphi,S)(u)$, $u \in [0,\bar{U}]$.
\end{Lemma}
{\bf Proof.} We have 
\begin{align*}\widehat{S}(t, z\vert w)&=\widehat{p}_{zw} \widehat{\widetilde{S}}(t|z,w)\\
&= \widehat{p}_{zw}\int \widehat{S}_{KM}(t-s\epsilon\vert z, w)K(s)ds\\
&=-\widehat{p}_{zw}\sum_{i=1}^{N_{\delta}}\widehat{S}_{(i)} \Big[\bar{K}\Big(\frac{t-Y_{(i+1)}}{\epsilon}\Big)-\bar{K}\Big(\frac{t-Y_{(i)}}{\epsilon}\Big)\Big],
\end{align*}
where $Y_{(i)}$ is the $i^{th}$ order statistic of the sample of non-censored observations $\{Y_i|\ \delta_i=1\}_{i=1}^{n}$, $\widehat{S}_{(i)}=\widehat{S}_{KM}(Y_{(i)}\vert z,w) $, $\bar{K}(y)= \int_{0}^yK(s) ds$ and $ N_{\delta}=\sum_{i=1}^nI(\delta_i=1)$.
This implies that the mapping
$ t\mapsto \widehat{S}(t, z\vert w)$ is differentiable in $t$. Therefore, the mapping $\theta \mapsto A(\theta,\widehat{S})$, $\theta \in \mathcal{F}_Z$, is G\^ateaux differentiable in $\theta$. By Lemma \ref{finalcv}, $\widehat{S}$ converges in probability to $S$. The map $h \mapsto A(\varphi,h)$, $h \in \mathcal{F}_{Z,W}$, and the differentiation map, are continuous mappings, therefore, by the continuous mapping theorem, we obtain the result. The proof for the bootstrap is similar. \hfill $\Box$\\[-.3cm]

\begin{Lemma}\label{cvderiv2} 
Under Assumption (K) and if $\widehat{\varphi}\xrightarrow{P}\varphi$, then
$\Gamma(\widehat{\varphi},\widehat{S})(u)$, $u \in [0,\bar{U}]$ exists and converges in probability to $\Gamma(\varphi,S)$. Also, if $\widehat{\varphi}_b\xrightarrow{P}\varphi$, then $\Gamma(\widehat{\varphi}_b,\widehat{S}_b)(u)$, $u \in [0,\bar{U}]$ exists and converges in probability to $\Gamma(\varphi,S)$.
\end{Lemma}
{\bf Proof.} Similarly as in the proof of Lemma \ref{cvderiv}, $\theta \in \mathcal{F}_Z\mapsto A(\theta,\widehat{S})$  is G\^ateaux differentiable in $\theta$. By Lemma \ref{finalcv}, $\widehat{S}$ converges in probability to $S$. The map $A$ and the differentiation map are continuous mappings, therefore, by the continuous mapping theorem, we obtain the result. The proof for the bootstrap is similar. 
 \hfill $\Box$\\[-.3cm]
 
 \begin{Lemma}\label{derivons} 
 Under Assumption (K),  $A(\widehat{\varphi},\widehat{S})-A(\varphi,\widehat{S})-\Gamma(\varphi,\widehat{S})(\widehat{\varphi}-\varphi)=O_P(\| \widehat{\varphi}-\varphi \|^2).$ Also, $A(\widehat{\varphi}_b,\widehat{S}_b)-A(\widehat{\varphi},\widehat{S}_b)-\Gamma(\widehat{\varphi},\widehat{S}_b)(\widehat{\varphi}_b-\widehat{\varphi})=O_P(\| \widehat{\varphi}_b-\widehat{\varphi} \|^2).$
\end{Lemma}
{\bf Proof.} Similarly as in the proof of Lemma \ref{cvderiv}, $\theta \mapsto A(\theta,\varphi)$, $\theta \in \mathcal{F}_Z$,  is twice G\^ateaux differentiable in $\theta$. By Lemma \ref{finalcv}, $\widehat{S}$ converges in probability to $S$. The map $A$ and the differentiation map are continuous mappings, therefore, by the continuous mapping theorem, we obtain that the second differential of $A$ with respect to its first argument in $(\varphi,\widehat{S})$  converges in probability to the second differential of $A$ with respect to its first argument in $(\varphi,S)$ which is bounded by Assumption (K)(i). This implies the result by a Taylor expansion. The proof for the bootstrap is similar.\hfill $\Box$\\[-.3cm]
 
 \begin{Lemma}\label{boot} 
 Under the assumptions of Lemma \ref{K}, the bootstrap approximation of Section 4.4 works in the sense that $\sqrt{n}(\widehat{\varphi}_b-\widehat{\varphi})$ converges weakly to the same mean zero Gaussian process as $\sqrt{n}(\widehat{\varphi}-\varphi)$. 
\end{Lemma}
{\bf Proof.} By lemmas \ref{FAN} and \ref{cvS}, $\sqrt{n}(\widehat{S}-S)$ and $\sqrt{n}(\widehat{S}_b-\widehat{S})$ converge to the same mean zero Gaussian process. Therefore, by the functional Delta-method, as $A$ is  a linear mapping, $\sqrt{n}A(\varphi,\widehat{S})$ and  $\sqrt{n}(A(\varphi,\widehat{S}_b)-A(\varphi,\widehat{S}))$ converge weakly to the same mean zero Gaussian process with variance operator $\Omega$.  Using Theorem \ref{C}, we obtain that $\|\widehat{\varphi}_b-\varphi\|=O_P(n^{-1/2})$ ((C)(ii) holds with $r_n=n^{-1/2}$ in the resample $b$ because of Lemma \ref{finalcv}) and $\|\widehat{\varphi}-\varphi\|=O_P(n^{-1/2})$. This implies that $\|\widehat{\varphi}_b-\widehat{\varphi}\|=O_P(n^{-1/2})$. The first-order conditions of the optimisation program \eqref{prog} in resample $b$ are 
$$\Gamma(\widehat{\varphi}_b,\widehat{S}_b)^\top V A(\widehat{\varphi}_b,\widehat{S}_b)=0,$$
which leads to 
\begin{align*}
& \Gamma(\widehat{\varphi}_b,\widehat{S}_b)^\top VA(\widehat{\varphi},\widehat{S}_b) 
+ \Gamma(\widehat{\varphi}_b,\widehat{S}_b)^\top V\Gamma(\widehat{\varphi},\widehat{S}_b)(\widehat{\varphi}_b-\widehat{\varphi})
+ \Gamma(\widehat{\varphi}_b,\widehat{S}_b)^\top V R=0,
\end{align*}
where $R= A(\widehat{\varphi}_b,\widehat{S}_b)-A(\widehat{\varphi},\widehat{S}_b)-\Gamma(\widehat{\varphi},\widehat{S}_b)(\widehat{\varphi}_b-\widehat{\varphi})$.
Then, 
\begin{align*}
\sqrt{n} (\widehat{\varphi}_b-\widehat{\varphi})
= & -[ \Gamma(\widehat{\varphi}_b,\widehat{S}_b)^\top V\Gamma(\widehat{\varphi},\widehat{S}_b)]^{-1}\Gamma(\widehat{\varphi}_b,\widehat{S}_b)^\top V\sqrt{n}A(\widehat{\varphi},\widehat{S}_b)\\
& - [ \Gamma(\widehat{\varphi}_b,\widehat{S}_b)^\top V\Gamma(\widehat{\varphi},\widehat{S}_b)]^{-1} \Gamma(\widehat{\varphi}_b,\widehat{S}_b)^\top V\sqrt{n} R\\
= & -[ \Gamma(\widehat{\varphi}_b,\widehat{S}_b)^\top V\Gamma(\widehat{\varphi},\widehat{S}_b)]^{-1}\Gamma(\widehat{\varphi}_b,\widehat{S}_b)^\top V\sqrt{n}(A(\widehat{\varphi},\widehat{S}_b)-A(\widehat{\varphi},\widehat{S}))\\
& -[ \Gamma(\widehat{\varphi}_b,\widehat{S}_b)^\top V\Gamma(\widehat{\varphi},\widehat{S}_b)]^{-1}\Gamma(\widehat{\varphi}_b,\widehat{S}_b)^\top V\sqrt{n}A(\widehat{\varphi},\widehat{S})\\
& - [ \Gamma(\widehat{\varphi}_b,\widehat{S}_b)^\top V\Gamma(\widehat{\varphi},\widehat{S}_b)]^{-1} \Gamma(\widehat{\varphi}_b,\widehat{S}_b)^\top V\sqrt{n} R.
\end{align*}
Now, by Lemma \ref{derivons} and the fact that $\|\widehat{\varphi}_b-\widehat{\varphi}\|=O_P(n^{-1/2})$, we have $\sqrt{n}R=O_P(n^{1/2}n^{-1})=o_P(1)$, and hence 
$$[ \Gamma(\widehat{\varphi}_b,\widehat{S}_b)^\top V\Gamma(\widehat{\varphi},\widehat{S}_b)]^{-1} \Gamma(\widehat{\varphi}_b,\widehat{S}_b)^\top V \sqrt{n} R=o_P(1),$$ by Lemmas \ref{cvderiv} and \ref{cvderiv2} and Assumption (V). Moreover, as $$\Gamma(\widehat{\varphi},\widehat{S})^\top V A(\widehat{\varphi},\widehat{S})=0,$$
we have $$ \Gamma(\widehat{\varphi}_b,\widehat{S}_b)^\top V\sqrt{n}A(\widehat{\varphi},\widehat{S})=(\Gamma(\widehat{\varphi}_b,\widehat{S}_b)-\Gamma(\widehat{\varphi},\widehat{S}))^\top V\sqrt{n}A(\widehat{\varphi},\widehat{S}).$$
Now, by Lemma \ref{derivons} and the fact that $\|\widehat{\varphi}-\varphi\|=O_P(n^{-1/2})$, we have that $\sqrt{n}A(\widehat{\varphi},\widehat{S})=\sqrt{n}\Gamma(\varphi,\widehat{S})(\widehat{\varphi}-\varphi)+O_P(1)=O_P(1)$ where the last equality is because $\Gamma(\varphi,\widehat{S})=O_P(1)$ by Lemma \ref{cvderiv}. By Lemmas \ref{cvderiv} and \ref{cvderiv2} and Assumption (V), we obtain that $$[ \Gamma(\widehat{\varphi}_b,\widehat{S}_b)^\top V\Gamma(\widehat{\varphi},\widehat{S}_b)]^{-1}\Gamma(\widehat{\varphi}_b,\widehat{S}_b)^\top V\sqrt{n}A(\widehat{\varphi},\widehat{S})=o_P(1).$$
Therefore, we have 
$$\sqrt{n} (\widehat{\varphi}_b-\widehat{\varphi})=-[ \Gamma(\widehat{\varphi}_b,\widehat{S}_b)^\top V\Gamma(\widehat{\varphi},\widehat{S}_b)]^{-1}\Gamma(\widehat{\varphi}_b,\widehat{S}_b)^\top V\sqrt{n}(A(\widehat{\varphi},\widehat{S}_b)-A(\widehat{\varphi},\widehat{S}))+o_P(1).$$
Let us denote by $\frac{\partial\widehat{S}}{\partial t} $ (resp. $\frac{\partial\widehat{S}_b}{\partial t}$), the derivative of $\widehat{S}$ (resp. $\widehat{S}_b$) in its first argument. Next, by Lemma \ref{derivons}, we have \begin{align*}& \sqrt{n}(A(\widehat{\varphi},\widehat{S}_b)-A(\widehat{\varphi},\widehat{S}))-\sqrt{n}(A(\varphi,\widehat{S}_b)-A(\varphi,\widehat{S}))\\
&=\sqrt{n}\Gamma(\varphi,\widehat{S}_b)(\widehat{\varphi}-\varphi)- \sqrt{n}\Gamma(\varphi,\widehat{S})(\widehat{\varphi}-\varphi)+o_P(1)\\\
&=\sqrt{n}(\sum_{\ell=1}^L \frac{\partial\widehat{S}_b}{\partial t}(\varphi(z_{\ell},\cdot),z_\ell\vert w_k)(\widehat{\varphi}(z_{\ell},\cdot)-\varphi(z_{\ell},\cdot)))_{k=1}^K\\
&\quad - \sqrt{n}(\sum_{\ell=1}^L \frac{\partial\widehat{S}}{\partial t}(\varphi(z_{\ell},\cdot),z_\ell\vert w_k)(\widehat{\varphi}(z_{\ell},\cdot)-\varphi(z_{\ell},\cdot)))_{k=1}^K+o_P(1)\\
&=\sqrt{n}(\sum_{\ell=1}^L (\frac{\partial\widehat{S}_b}{\partial t}(\varphi(z_{\ell},\cdot),z_\ell\vert w_k)-\frac{\partial\widehat{S}}{\partial t}(\varphi(z_{\ell},\cdot),z_\ell\vert w_k))(\widehat{\varphi}(z_{\ell},\cdot)-\varphi(z_{\ell},\cdot)))_{k=1}^K+o_P(1).
\end{align*}
Now, by Lemma \ref{FAN}, $\widehat{S}_b$ converges to $\widehat{S}$ in probability. As the differentiation is a continuous mapping, we obtain that, for $z=z_1,\dots,z_L, w=w_1,\dots,w_K$, 
$$ \Big\|\frac{\partial\widehat{S}_b}{\partial t}(\varphi(z,\cdot),z\vert w)-\frac{\partial\widehat{S}}{\partial t}(\varphi(z,\cdot),z\vert w)\ \Big\|=o_P(1).$$
This implies that 
$$\sqrt{n}(A(\widehat{\varphi},\widehat{S}_b)-A(\widehat{\varphi},\widehat{S}))-\sqrt{n}(A(\varphi,\widehat{S}_b)-A(\varphi,\widehat{S}))=o_P(1).$$
Therefore,  $\sqrt{n}(A(\widehat{\varphi},\widehat{S}_b)-A(\widehat{\varphi},\widehat{S}))$ converge weakly to the mean zero Gaussian process with variance operator $\Omega$. We conclude that $\sqrt{n}(\widehat{\varphi}_b-\widehat{\varphi})$ converges weakly to a mean zero Gaussian process with variance operator
$[\Sigma^\top V\Sigma]^{-1} \Sigma^\top V \Omega V\Sigma[\Sigma^\top V\Sigma]^{-1} $ using Slutsky's theorem and Lemmas \ref{cvderiv} and \ref{cvderiv2}. By Theorem \ref{AS} (which we can apply because of the previous lemmas of this section), $\sqrt{n}(\widehat{\varphi}-\varphi)$ converges to the same Gaussian process, which yields the result.
\hfill $\Box$\\

\section{Appendix C: Primitive conditions for Assumption (C) (ii)}
\label{app.C}
\begin{Lemma}
\label{sufficient conditions}
Assume that $S(\cdot,z \vert w)$ is twice differentiable and has bounded first and second derivatives on $[0,\bar{T}]$, then 
Assumption (C) (ii) holds.
\end{Lemma}
{\bf Proof.} Because $S(\cdot,z \vert w)$ is twice differentiable on $[0,\bar{T}]$, $A(\cdot,S)$ is twice differentiable in its first argument. Also, its second derivative is bounded in $\mathcal{F}_Z^{\bar{U},\bar{T}}$ because $S(\cdot,z \vert w)$ has bounded second derivative on $[0,\bar{T}]$. By a Taylor expansion, we obtain that there exists a constant $H_1>0$ such that, for $\mathcal{F}_Z^{\bar{U},\bar{T}}$, we have 
$$\|A(\theta,S)-A(\varphi,S)- \Gamma(\varphi,S) (\theta-\varphi)\|\le H_1\|\theta-\varphi\|^2 .$$
As $A(\varphi,S)=0$, we obtain 
$$\| \Gamma(\varphi,S) (\theta-\varphi)\|\le \| A(\theta,S) \|+H_1\|\theta-\varphi\|^2 .$$
Because $S(\cdot,z \vert w)$ has bounded derivative on $[0,\bar{T}]$, $\Gamma(\varphi,S)$ is bounded on $[0,\bar{U}]$ and, therefore, there exists as constant $H_2>0$ such that 
$$H_2\| \theta-\varphi\| -H_1\|\theta-\varphi\|^2\le \| A(\theta,S) \| .$$
By choosing $\nu>0$ small enough, we have that, if $\| \theta-\varphi\|\le \nu $,$$\frac{H_2}{2}\| \theta-\varphi\| \le \| A(\theta,S) \| .$$
\vspace*{-1.5cm}
\hfill $\Box$\\

\end{document}